\documentclass[12pt,draft,leqno]{article}
\usepackage{amssymb, eucal, latexsym}

\newtheorem{theorem}{Theorem}[section]
\newtheorem{lm}[theorem]{Lemma}
\newtheorem{exa}[theorem]{Example}

\newtheorem{cor}[theorem]{Corollary}
\newtheorem{pro}[theorem]{Proposition}
\newtheorem{defi}[theorem]{Definition}

\newtheorem{nota}[theorem]{Notation}

\newtheorem{rem}[theorem]{Remark}

\newtheorem{fact}[theorem]{Fact}

\newtheorem{nist}[theorem]{}

\def\p{\varphi}
\def\a{\alpha}

\def\d{\delta}
\def\ep{\varepsilon}
\def\g{\gamma}
\def\GA{\Gamma}
\def\DE{\Delta}

\def\l{\lambda}

\def\t{\theta}

\def\s{\sigma}

\def\lra{\longrightarrow}

\def\sbe{\subseteq}

\def\stm{\setminus}
\def\ems{\emptyset}
\def\nes{\neq\emptyset}

\def\cuk{\,\check{}\,}

\def\ti{\tilde}
\def\wt{\widetilde}

\def\unl{\underline}
\def\ts{\tilde{\sigma}}

\def\ex{\exists}
\def\fa{\forall}

\def\we{\wedge}
\def\bw{\bigwedge}
\def\bv{\bigvee}

\def\ap{^\prime}
\def\inv{^{-1}}
\def\st{\ |\ }

\def\llx{\ll_{\rho}}
\def\lle{\ll_{\eta}}

\def\nin{\not\in}

\def\card #1{\vert #1 \vert}

\def\ion{i=1,\ldots,n}
\def\jom{j=1,\ldots,m}

\def\oon{\{1,\ldots,n\}}

\def\BB{{\cal B}}
\def\CC{{\cal C}}

\def\HC{{\bf HC}}
\def\DHC{{\bf DHC}}

\def\HLC{{\bf HLC}}
\def\DHLC{{\bf DHLC}}

\def\Bool{{\bf Bool}}

\def\A{\mbox{{\boldmath $A$}}}
\def\B{\mbox{{\boldmath $B$}}}

\def\T{\mbox{{\boldmath $T$}}}

\def\2{\mbox{{\bf 2}}}
\def\3{\mbox{{\bf 3}}}

\def\int{\mbox{{\rm int}}}

\def\cl{\mbox{{\rm cl}}}

\def\doc{\hspace{-1cm}{\em Proof.}~~}
\def\sq{\hspace*{\fill} \hbox{\vrule\vbox{\hrule\phantom{o}\hrule}\vrule}}
\def\sqs{\sq \vspace{2mm}}

\def\NNNN{{\rm I}\!{\rm N}}
\def\BBBB{{\rm I}\!{\rm B}}

\def\eset{\emptyset}   %empty set%
\def\neset{\neq\emptyset} %%not=emptyset%%

\title{{\LARGE\bf
A Whiteheadian-type description of Euclidean spaces, spheres, tori and Tychonoff cubes}\\
\vspace{0.35cm}
{\large\bf Georgi D. Dimov}\thanks{This paper was supported by the
project no. DID 02/32/2009  $``$Theories of the space and time: algebraic, topological and logical approaches" of the Bulgarian Ministry of Education and Science.
}\\
\vspace{0.25cm}
 {\footnotesize Dept. of Math. and
Informatics, Sofia University,  5 J. Bourchier Blvd., 1164 Sofia,
Bulgaria}
}

\author{}

\date{}

\begin{document}
\maketitle
\begin{abstract}
In the beginning of the 20th century, A. N. Whitehead \cite{W1,W} and T. de Laguna \cite{dL} proposed
a new theory of space, known as {\em region-based theory of space}.
They did not present  their ideas in a detailed mathematical form.
 In 1997,
 P. Roeper \cite{R} has shown that the locally compact
Hausdorff spaces correspond bijectively (up to homeomorphism and isomorphism) to some algebraical objects
which represent correctly  Whitehead's idea of {\em regions}\/ and {\em contact relation},
generalizing in this way a previous analogous result of de Vries \cite{dV} concerning  compact Hausdorff spaces    (note that even a duality  for the category of compact Hausdorff spaces and continuous maps was constructed by de Vries \cite{dV}). Recently,  a duality for the category of locally compact Hausdorff spaces and continuous maps, based on  Roeper's results, was obtained in \cite{DI} (it extends de Vries' duality mentioned above). In this paper,
 using the  dualities obtained in \cite{dV,DI}, we construct directly (i.e. without the help of the corresponding topological spaces) the dual  objects of  Euclidean spaces, spheres, tori and Tychonoff cubes; these algebraical objects completely characterize the mentioned topological spaces. Thus,  a mathematical realization of the original philosophical ideas of Whitehead \cite{W1,W} and de Laguna \cite{dL} about Euclidean spaces is obtained.
\end{abstract}



\footnotetext[3]{{\footnotesize {\em E-mail address:}
gdimov@fmi.uni-sofia.bg}}

%\vspace{1.5cm}

\baselineskip = \normalbaselineskip

\section{Introduction}

%\centerline{{\bf 1. Introduction}}

\bigskip

The region-based theory of space is a kind of
point-free geometry and can be considered as  an alternative to
the well known Euclidean point-based theory of space. Its main
idea goes back to Whitehead \cite{W} (see also \cite{W1}) and de
Laguna \cite{dL} and is based on a certain criticism of the
Euclidean approach to the geometry, where the points (as well as
straight lines and planes) are taken as the basic primitive
notions. A. N. Whitehead and T. de Laguna noticed that points, lines and planes are quite
abstract entities which have not a separate existence in reality
and proposed to put the theory of space  on the base of some more
realistic spatial entities. In Whitehead \cite{W}, the notion of a
{\em region}\/ is taken as a primitive notion: it is an abstract analog of
a spatial body; also some natural relations between regions are
regarded. In \cite{W1}, Whitehead considered  some mereological
relations like ``part-of", ``overlap" and some others, while in \cite{W} he
adopted  from de Laguna \cite{dL} the relation of ``{\em contact}\/"
(``connectedness" in Whitehead's original terminology) as the only
primitive relation between regions except the relation ``part-of". The {\em regular closed}\/ (or,
equivalently, {\em regular open}\/) subsets of a topological space $X$
are usually considered as a standard model of the regions in the
point-based approach, and the {\em standard contact relation} $\rho_X$ between regular closed subsets of $X$
is defined (again in the point-based approach) as follows:  $F\rho_X G\Leftrightarrow F\cap G\nes$.

 Let us note that neither Whitehead nor de
Laguna presented  their ideas in a detailed mathematical form. This was done by some other mathematicians and mathematically oriented
philosophers who presented various versions of region-based theory of
space at different levels of abstraction. Here we can mention
Tarski \cite{Tarski}, who rebuilt Euclidean geometry as an
extension of mereology with the primitive notion of a {\em ball}.
Remarkable is also Grzegorczyk's paper \cite{Grzegorczyk}. Models of
Grzegorczyk's  theory  are complete Boolean algebras of regular
closed sets of certain topological spaces equipped with the
relation of separation which in fact is the complement of
Whitehead's contact relation. On the same line of abstraction is
also the point-free topology \cite{J}. Survey papers describing
various aspects and historical remarks on region-based theory of
space are \cite{Gerla,BD,Vak,Pratt}.

Let us mention that Whitehead's ideas about region-based theory of
space  flourished and in a sense  were reinvented and applied in
some areas of computer science: Qualitative Spatial Reasoning
(QSR), knowledge representation, geographical information systems,
formal ontologies in information systems, image processing,
natural language semantics etc. The reason is that the language of
region-based theory of space allows the researches to obtain a more simple  description of
some qualitative spatial features and properties of space bodies.
Survey papers concerning various applications are \cite{ch01,CR}
(see also the special issues of $``$Fundamenta Informaticae" \cite{FI} and
 $``$Journal of Applied Non-classical Logics" \cite{JANCL}). One of
the most popular among the community of QSR-researchers is the
system of Region Connection Calculus (RCC) introduced by Randell,
Cui and Cohn \cite{Randell}. RCC attracted quite intensive
research in the field of region-based theory of space, both on its
applied and mathematical aspects. For instance it was unknown for
some time which topological models correspond adequately to RCC;
this fact stimulated the investigations of a topological
representation theory of RCC and RCC-like systems (see
\cite{DW,DV}). Another impact of region-based theory of space is
that it stimulated the appearance of a new area in logic, namely ``Spatial Logics"
 \cite{A}, called sometimes ``Logics of Space".

The ideas of de Laguna and Whitehead lead naturally to the
following general programme (or {\em general region-based theory of space}):

\begin{itemize}
 \item for every topological space $X$ belonging to some class $\mathcal{C}$ of topological spaces, define in topological terms:

 (a) a family $\mathcal{R}(X)$ of subsets of $X$
 that will serve as models of   Whitehead's 
$``$regions" (and call the elements of the family $\mathcal{R}(X)$ {\em regions of $X$});

(b) a relation $\rho_X$ on $\mathcal{R}(X)$  that will serve as a model of  Whitehead's relation of
$``$contact" (and call the relation $\rho_X$ a {\em contact relation on} $\mathcal{R}(X)$);
 \item choose some (algebraic) structure which is inherent
to the families $\mathcal{R}(X)$ and contact relations $\rho_X$, for $X\in\mathcal{C}$, fix some kind
of morphisms between the obtained (algebraic) objects and build in
this way a category $\A$;
\item find a subcategory $\T$ of the
category of topological spaces and continuous maps which is
equivalent or dually equivalent to the category $\A$ trough a
(contravariant) functor that assigns to each object $X$ of $\T$
the chosen (algebraic) structure of the family of all regions of
$X$.
\end{itemize}

 If all of this is done then, in particular, the chosen
(algebraic) structure of the regions of any object $X$ of $\T$ is
sufficient for recovering completely (of course, up to
homeomorphism) the whole space $X$. Hence, in this way, a
$``$region-based theory" of the objects and morphisms of the
category $\T$ is obtained.

Of course, during the realization of this programme, one can find
the category $\A$ starting with the category $\T$, if the later is
the desired one.

 The M. Stone \cite{ST} duality between
 the category of Boolean algebras and their homomorphisms and
the category of compact zero-dimensi\-onal Hausdorff spaces and
continuous maps can be regarded as a first realization of this programme, although M. Stone came to his results guided by ideas which are completely different from those of Whitehead and de Laguna. In M. Stone's theory, the clopen (= closed and open) subsets of a topological
space serve as models of the regions; here, however,  the  contact relation $\rho$ is hidden, because it can be defined by the
 Boolean operations (indeed, we have that $a\rho b\iff a\we b\neq 0$). The {\em localic duality} (see, e.g., \cite[Corollary II.1.7]{J}) between the category of spatial frames and functions preserving finite meets and arbitrary joins and the category of sober spaces and continuous maps can also be regarded as a realization of the ideas of the
 general region-based theory of space:  in it the open subsets of a topological space serve as models of the regions and, as above, the  contact relation $\rho$ between the regions is hidden because it can be recovered by the lattice operations (indeed, we have that $a\rho b\iff a\we b\neq 0$).
 The de Vries   duality \cite{dV} for the category $\HC$ of compact Hausdorff spaces and continuous maps
is the first  realization of the ideas of the general region-based theory of space in their full generality and strength
 (and again, as it seems, de Vries was unaware of
the papers  \cite{dL} and  \cite{W}): the models of the regions in de
Vries' theory are the regular closed sets and, in contrast to
the case of the Stone duality and localic duality, the  contact
relation between regions, which is in the basis of de Vries' duality theorem, cannot be
derived from the Boolean structure on the regions.
(Note that in \cite{dV}, instead of the Boolean algebra $RC(X)$ of regular closed sets, the Boolean algebra $RO(X)$ of regular open sets was regarded ($RO(X)$ and $RC(X)$ are isomorphic); also, instead of the relation $\rho_X$ on the set $RC(X)$ which was described above (let us recall it: $F\rho_X G\iff  F\cap G\nes$), de Vries used in  \cite{dV} the so-called {\em $``$compingent relation"}\/ between regular open sets whose counterpart for $RC(X)$ is
the relation $\ll_X$, defined by $F\ll_X G\iff F\sbe
\int(G)$, for $F,G\in RC(X)$;  the relations $\rho_X$ and $\ll_X$ are inter-definable.)
It is natural to try to extend de Vries' Duality Theorem to the
category $\HLC$ of locally compact Hausdorff spaces and continuous
maps. An important step in this direction was done by P. Roeper
\cite{R}. Being guided by the ideas of de Laguna \cite{dL} and
Whitehead \cite{W}, he  proved that there is a
bijective correspondence between all (up to homeomorphism) locally
compact Hausdorff spaces and all (up to isomorphism) algebras of some sort called by him {\em $``$region-based topologies"}\/ (we call them {\em complete
LC-algebras}). The notion of a complete LC-algebra, introduced by Roeper \cite{R}, is an abstraction of the triples $(RC(X), \rho_X,CR(X))$, where $X$ is a locally
compact Hausdorff space and $CR(X)$ is the ideal of all  compact regular closed subsets of $X$. P. Roeper
\cite{R} showed that every complete LC-algebra can be realized  as a triple $(RC(X), \rho_X,CR(X))$, where $X$ is a uniquely (up to homeomorphism) determined locally compact Hausdorff
space.   In \cite{DI}, using Roeper's result, we obtained a duality
between the category $\HLC$ and the category $\DHLC$ of complete LC-algebras and appropriate morphisms between them; it is an extension of   de Vries'
duality mentioned above; the dual object of a locally compact Hausdorff
space $X$ is the triple $(RC(X), \rho_X,CR(X))$ which will be called {\em the Roeper triple of the space} $X$.  Let us note that the famous Gelfand duality \cite{G1,G2,GN,GS} also gives an algebraical description of (locally) compact Hausdorff spaces but it is not in the spirit of the ideas of Whitehead and de Laguna.

A description of the dual object of the real line under the localic duality (i.e., a description of the frame (or locale) determined by the topology of the real line) without the help of the real line was given by Fourman and Hyland \cite{FH} (see, also, Grayson \cite{Gr} and Johnstone \cite[IV.1.1-IV.1.3]{J}), assuming the set of rationals as given. As we have seen above, the ideas of the localic duality are in the spirit of the ideas of the {\em general}\/ region-based theory of space but, nevertheless, they are far from the well-known and commonly accepted interpretations  of the {\em original}\/ philosophical ideas of Whitehead \cite{W1,W} and de Laguna \cite{dL} given in \cite{Grzegorczyk} and \cite{R} (see also \cite{Randell}).

In this paper we construct directly the dual objects of  Euclidean spaces, spheres, tori and Tychonoff cubes under the dualities obtained in \cite{dV,DI}, i.e. we construct   the complete LC-algebras  isomorphic to the Roeper triples (see \cite{R}) of these spaces   without the help of the corresponding spaces, assuming the set of natural numbers as given. For doing this, we first obtain some direct descriptions of  the $\DHLC$-sums of complete LC-algebras and the $\DHC$-sums of complete NC-algebras (where $\DHC$ is the de Vries category dual to the category $\HC$, and  the objects of the category $\DHC$  are the complete NC-algebras) using the dualities obtained in \cite{dV} and \cite{DI}.  Let us note explicitly that, as it follows from the results of de Vries \cite{dV} and Roeper \cite{R}, the Euclidean spaces, spheres, tori and Tychonoff cubes can be completely reconstructed as topological spaces from the algebraical objects which we describe in this paper.  Therefore, our results can be regarded as a mathematical realization of the original philosophical ideas of Whitehead \cite{W1,W} and de Laguna \cite{dL} about Euclidean spaces; this realization is in accordance with the Grzegorczyk's \cite{Grzegorczyk} and Roeper's \cite{R} mathematical interpretations of these ideas.

We now fix the notation.

If $\CC$ denotes a category, we write $X\in \card\CC$ if $X$ is an
object of $\CC$, and $f\in \CC(X,Y)$ if $f$ is a morphism of $\CC$
with domain $X$ and codomain $Y$.

All lattices are with top (= unit) and bottom (= zero) elements,
denoted respectively by 1 and 0. We do not require the elements
$0$ and $1$ to be distinct.

 If $(X,\tau)$ is a topological space and $M$ is a subset of $X$, we
denote by $\cl_{(X,\tau)}(M)$ (or simply by $\cl(M)$ or
$\cl_X(M)$) the closure of $M$ in $(X,\tau)$ and by
$\int_{(X,\tau)}(M)$ (or briefly by $\int(M)$ or $\int_X(M)$) the
interior of $M$ in $(X,\tau)$. The Alexandroff compactification of
a locally compact Hausdorff non-compact space $X$ will be denoted
by $\a X$.
The positive natural numbers are denoted by  ${\NNNN^+}$,
  the real line (with its natural topology) -- by $\mathbb{R}$,  the $n$-dimensional sphere
  (with its natural topology) -- by $\mathbb{S}^n$ (here $n\in{\NNNN^+}$).

\section{Preliminaries}
%%%%%%%%%%%%%%%%%%%%%%%%%%%%%%%%%%%%%%%%%%%%%%%%%%%%%%%%%%%%
%%%% Section 1. %%%%%%%%%%%%%%%%%%%%%%%
%%%%%%%%%%%%%%%%%%%%%%%%%%%%%%%%%%%%%%%%%%%%%%%%%%%%%%%%%%%%
%

\begin{defi}\label{conalg}
\rm
An algebraic system $(B,0,1,\vee,\we, {}^*, C)$ is called a {\it
contact Boolean algebra}\/ or, briefly, {\it contact algebra}
(abbreviated as CA or C-algebra) (\cite{DV})
 if the system
$(B,0,1,\vee,\we, {}^*)$ is a Boolean algebra (where the operation
$``$complement" is denoted by $``\ {}^*\ $")
  and $C$
is a binary relation on $B$, satisfying the following axioms:

\smallskip

\noindent (C1) If $a\not= 0$ then $aCa$;\\
(C2) If $aCb$ then $a\not=0$ and $b\not=0$;\\
(C3) $aCb$ implies $bCa$;\\
(C4) $aC(b\vee c)$ iff $aCb$ or $aCc$.

\smallskip

\noindent We shall simply write $(B,C)$ for a contact algebra. The
relation $C$  is called a {\em  contact relation}. When $B$ is a
complete Boolean algebra, we will say that $(B,C)$ is a {\em
complete contact Boolean algebra}\/ or, briefly, {\em complete
contact algebra} (abbreviated as CCA or CC-algebra). If $a\in B$ and $D\sbe B$,
we will write $``aCD$" for $``(\fa d\in D)(aCd)$".

We will say that two C-algebras $(B_1,C_1)$ and $(B_2,C_2)$ are  {\em
CA-isomorphic} iff there exists a Boolean isomorphism $\p:B_1\lra
B_2$ such that, for each $a,b\in B_1$, $aC_1 b$ iff $\p(a)C_2
\p(b)$. Note that in this paper, by a $``$Boolean isomorphism" we
understand an isomorphism in the category $\Bool$ of Boolean algebras and Boolean homomorphisms.

\smallskip

A contact algebra $(B,C)$ is called a {\it  normal contact Boolean
algebra}\/ or, briefly, {\it  normal contact algebra} (abbreviated
as NCA or NC-algebra) (\cite{dV,F}) if it satisfies the following axioms which are
very similar to the Efremovi\v c \cite{EF} axioms of proximity spaces
(we
will write $``-C$" for $``not\ C$"):

\smallskip

\noindent (C5) If $a(-C)b$ then $a(-C)c$ and $b(-C)c^*$ for some $c\in B$;\\
(C6) If $a\not= 1$ then there exists $b\not= 0$ such that
$b(-C)a$.

\smallskip

\noindent A normal CA is called a {\em complete normal contact
Boolean algebra}\/ or, briefly, {\em complete normal contact
algebra} (abbreviated as CNCA or CNC-algebra) if it is a CCA. The notion of a
normal contact algebra was introduced by Fedorchuk \cite{F} under
the name {\em Boolean $\d$-algebra}\/ as an equivalent expression
of the notion of a {\em compingent Boolean algebra}\/ of de Vries (see its definition below). We call
such algebras $``$normal contact algebras" because they form a
subclass of the class of contact algebras and naturally arise in
normal Hausdorff spaces.

Note that if $0\neq 1$ then the axiom (C2) follows from the axioms
(C6) and (C4).

For any CA $(B,C)$, we define a binary relation  $``\ll_C $"  on
$B$ (called {\em non-tangential inclusion})  by $``\ a \ll_C b
\leftrightarrow a(-C)b^*\ $". Sometimes we will write simply
$``\ll$" instead of $``\ll_C$".
\end{defi}

The relations $C$ and $\ll$ are inter-definable. For example,
normal contact algebras could be equivalently defined (and exactly
in this way they were introduced (under the name of {\em
compingent Boolean algebras}) by de Vries in \cite{dV}) as a pair
of a Boolean algebra $B=(B,0,1,\vee,\we,{}^*)$ and a binary
relation $\ll$ on $B$ subject to the following axioms:

\smallskip

\noindent ($\ll$1) $a\ll b$ implies $a\leq b$;\\
($\ll$2) $0\ll 0$;\\
($\ll$3) $a\leq b\ll c\leq d$ implies $a\ll d$;\\
($\ll$4) $a\ll c$ and $b\ll c$ implies $a\vee b\ll c$;\\
($\ll$5) If  $a\ll c$ then $a\ll b\ll c$  for some $b\in B$;\\
($\ll$6) If $a\neq 0$ then there exists $b\neq 0$ such that $b\ll
a$;\\
($\ll$7) $a\ll b$ implies $b^*\ll a^*$.

\smallskip

Note that if $0\neq 1$ then the axiom ($\ll$2) follows from the
axioms ($\ll$3), ($\ll$4), ($\ll$6) and ($\ll$7).

\smallskip

Obviously, contact algebras could be equivalently defined as a
pair of a Boolean algebra $B$ and a binary relation $\ll$ on $B$
subject to the  axioms ($\ll$1)-($\ll$4) and ($\ll$7).

\smallskip

It is easy to see that axiom (C5) (resp., (C6)) can be stated
equivalently in the form of ($\ll$5) (resp., ($\ll$6)).

\begin{exa}\label{rct}
\rm Recall that a subset $F$ of a topological space $(X,\tau)$ is
called {\em regular closed}\/ if $F=\cl(\int (F))$. Clearly, $F$
is regular closed iff it is the closure of an open set.

For any topological space $(X,\tau)$, the collection $RC(X,\tau)$
(we will often write simply $RC(X)$) of all regular closed subsets
of $(X,\tau)$ becomes a complete Boolean algebra
$(RC(X,\tau),0,1,\we,\vee,{}^*)$ under the following operations:
$$ 1 = X,  0 = \emptyset, F^* = \cl(X\stm F), F\vee G=F\cup G,
F\we G =\cl(\int(F\cap G)).
$$
The infinite operations are given by the  formulae:
$\bigvee\{F_\g\st \g\in\GA\}=\cl(\bigcup\{F_\g\st
\g\in\GA\})(=\cl(\bigcup\{\int(F_\g)\st \g\in\GA\})),$ and
$\bigwedge\{F_\g\st \g\in\GA\}=\cl(\int(\bigcap\{F_\g\st
\g\in\GA\})).$

It is easy to see that setting $F \rho_{(X,\tau)} G$ iff $F\cap
G\not = \ems$, we define a contact relation $\rho_{(X,\tau)}$ on
$RC(X,\tau)$; it is called a {\em standard contact relation}. So,
$(RC(X,\tau),\rho_{(X,\tau)})$ is a CCA (it is called a {\em
standard contact algebra}). We will often write simply $\rho_X$
instead of $\rho_{(X,\tau)}$. Note that, for $F,G\in RC(X)$,
$F\ll_{\rho_X}G$ iff $F\sbe\int_X(G)$.

Clearly, if $(X,\tau)$ is a normal Hausdorff space then the
standard contact algebra $(RC(X,\tau),\rho_{(X,\tau)})$ is a
complete NCA.

A subset $U$ of $(X,\tau)$ such that $U=\int(\cl(U))$ is said to
be {\em regular open}. The set of all regular open subsets of
$(X,\tau)$ will be denoted by $RO(X,\tau)$ (or briefly, by
$RO(X)$).
\end{exa}

%%%%%%%%%%%%%%%%%%%%%%%%%%%%%% %%%%%%%%%%%%%%%%%%%%%%%%%%%%%%%%%%%

The following notion is a lattice-theoretical counterpart of
Leader's notion of a {\em local proximity} (\cite{LE}):

\begin{defi}\label{locono}{\rm (\cite{R})}
\rm An algebraic system $\underline {B}_{\, l}=(B,0,1,\vee,\we,
{}^*, \rho, \BBBB)$ is called a {\it local contact Boolean
algebra}\/ or, briefly, {\it local contact algebra} (abbreviated
as LCA or LC-algebra)   if $(B,0,1, \vee,\we, {}^*)$ is a Boolean algebra,
$\rho$ is a binary relation on $B$ such that $(B,\rho)$ is a CA,
and $\BBBB$ is an ideal (possibly non proper) of $B$, satisfying
the following axioms:

\smallskip

\noindent(BC1) If $a\in\BBBB$, $c\in B$ and $a\ll_\rho c$ then
$a\ll_\rho b\ll_\rho c$ for some $b\in\BBBB$  (see Definition
\ref{conalg} for
$``\ll_\rho$");\\
(BC2) If $a\rho b$ then there exists an element $c$ of $\BBBB$
such that
$a\rho (c\we b)$;\\
(BC3) If $a\neq 0$ then there exists  $b\in\BBBB\stm\{0\}$ such
that $b\ll_\rho a$.

\smallskip

We shall simply write  $(B, \rho,\BBBB)$ for a local contact
algebra.  We will say that the elements of $\BBBB$ are {\em
bounded} and the elements of $B\stm \BBBB$  are  {\em unbounded}.
When $B$ is a complete Boolean algebra,  the LCA $(B,\rho,\BBBB)$
is called a {\em complete local contact Boolean algebra}\/ or,
briefly, {\em complete local contact algebra} (abbreviated as
CLCA or CLC-algebra).

We will say that two local contact algebras $(B,\rho,\BBBB)$ and
$(B_1,\rho_1,\BBBB_1)$ are  {\em LCA-isomorphic} if there exists a
Boolean isomorphism $\p:B\lra B_1$ such that, for $a,b\in B$,
$a\rho b$ iff $\p(a)\rho_1 \p(b)$, and $\p(a)\in\BBBB_1$ iff
$a\in\BBBB$. A map $\p:(B,\rho,\BBBB)\lra(B_1,\rho_1,\BBBB_1)$ is
called an {\em LCA-embedding} if $\p:B\lra B_1$ is an injective
Boolean homomorphism (i.e. Boolean monomorphism) and, moreover,
for any $a,b\in B$, $a\rho b$ iff $\p(a)\rho_1 \p(b)$, and
$\p(a)\in\BBBB_1$ iff $a\in\BBBB$.
\end{defi}

\begin{rem}\label{conaln}
\rm Note that if $(B,\rho,\BBBB)$ is a local contact algebra and
$1\in\BBBB$ then $(B,\rho)$ is a normal contact algebra.
Conversely, any normal contact algebra $(B,C)$ can be regarded as
a local contact algebra of the form $(B,C,B)$.
\end{rem}

\begin{defi}\label{Alexprn}{\rm (\cite{VDDB})}
\rm Let $(B,\rho,\BBBB)$ be a local contact algebra. Define a
binary relation $``C_{\rho,\BBBB}$" on $B$ by
\begin{equation}\label{crho}
aC_{\rho,\BBBB} b\ \mbox{ iff }\ a\rho b\ \mbox{ or }\ a,b\not\in\BBBB.
\end{equation}
It is called the\/ {\em Alexandroff extension of}\/ $\rho$ {\em
relatively to the LCA} $(B,\rho,\BBBB)$ (or, when there is no
ambiguity, simply, the \/ {\em Alexandroff extension of}\/
$\rho$).
\end{defi}

The following lemma
%from \cite{VDDB}
is a lattice-theoretical
counterpart of a theorem from Leader's paper \cite{LE}.

\begin{lm}\label{Alexprn1}{\rm (\cite{VDDB})}
Let $(B,\rho,\BBBB)$ be a local contact algebra. Then
$(B,C_{\rho,\BBBB})$, where $C_{\rho,\BBBB}$ is the Alexandroff extension of
$\rho$, is a normal contact algebra.
\end{lm}

\begin{nota}\label{compregn}
\rm Let $(X,\tau)$ be a topological space. We denote by
$CR(X,\tau)$ the family of all compact regular closed subsets of
$(X,\tau)$. We will often write  $CR(X)$ instead of $CR(X,\tau)$.
\end{nota}

\begin{fact}\label{stanlocn}{\rm (\cite{R})}
Let $(X,\tau)$ be a locally compact Hausdorff space. Then
 the triple
$(RC(X,\tau),\rho_{(X,\tau)}, CR(X,\tau))$
 (see Example \ref{rct} for $\rho_{(X,\tau)}$)
  is a complete local contact algebra; it is called a
{\em standard local contact algebra}.
\end{fact}

The next  theorem was proved by Roeper \cite{R} (but its
particular case concerning compact Hausdorff spaces and NC-algebras was
proved by de Vries \cite{dV}).

\begin{theorem}\label{roeperl}{\rm (P. Roeper \cite{R}
for locally compact spaces and de Vries \cite{dV} for compact
spaces)}
  There exists a bijective correspondence $\Psi^t$ between the
class of all (up
to homeomorphism) locally compact Hausdorff spaces and the class of all (up to isomorphism) CLC-algebras; its
restriction to the class of all (up to homeomorphism) compact Hausdorff spaces  gives a
bijective correspondence between the later class and the class of
all (up to isomorphism) CNC-algebras.
\end{theorem}

Let us recall  the definition of the correspondence
$\Psi^t$ mentioned in the above theorem:
if $(X,\tau)$ is a
locally compact Hausdorff space then
\begin{equation}\label{psit1}
\Psi^t(X,\tau)=(RC(X,\tau),\rho_{(X,\tau)},CR(X,\tau))
\end{equation}
(see Fact \ref{stanlocn} and Notation \ref{compregn} for the
notation).

\begin{defi}\label{dval}{\rm (De Vries \cite{dV})}
\rm Let  $\HC$ be the category of all compact Hausdorff spaces and
all continuous maps between them.

Let $\DHC$ be the category whose objects are all complete NC-algebras and
whose morphisms are all functions $\p:(A,C)\lra (B,C\ap)$ between
the objects of $\DHC$ satisfying the conditions:

\smallskip

\noindent(DVAL1) $\p(0)=0$;\\
(DVAL2) $\p(a\we b)=\p(a)\we \p(b)$, for all $a,b\in A$;\\
(DVAL3) If $a, b\in A$ and $a\ll_C b$, then $(\p(a^*))^*\ll_{C\ap}
\p(b)$;\\
(DVAL4) $\p(a)=\bigvee\{\p(b)\st b\ll_{C} a\}$, for every $a\in
A$,

\medskip

{\noindent}and let the composition $``\diamond$" of two morphisms
$\p_1:(A_1,C_1)\lra (A_2,C_2)$ and $\p_2:(A_2,C_2)\lra (A_3,C_3)$
of $\DHC$ be defined by the formula
\begin{equation}\label{diamc}
\p_2\diamond\p_1 = (\p_2\circ\p_1)\cuk,
\end{equation}
 where, for every
function $\psi:(A,C)\lra (B,C\ap)$ between two objects of $\DHC$,
$\psi\cuk:(A,C)\lra (B,C\ap)$ is defined as follows:
\begin{equation}\label{cukfc}
\psi\cuk(a)=\bigvee\{\psi(b)\st b\ll_{C} a\},
\end{equation}
for every $a\in A$.
\end{defi}

De Vries \cite{dV} proved the following duality theorem:

\begin{theorem}\label{dvth}{\rm (\cite{dV})}
The categories $\HC$ and $\DHC$ are dually equivalent.
\end{theorem}

In \cite{DI},  an extension of de Vries' Duality Theorem to the category of locally compact
 Hausdorff spaces and continuous maps was obtained. Let us recall its formulation.

\begin{defi}\label{dhc}{\rm (\cite{DI})}
\rm Let  $\HLC$ be the category of all locally compact Hausdorff
spaces and all continuous maps between them.

Let $\DHLC$ be the category whose objects are all complete LC-algebras
and whose morphisms are all functions $\p:(A,\rho,\BBBB)\lra
(B,\eta,\BBBB\ap)$ between the objects of $\DHLC$ satisfying
conditions

\smallskip

\noindent(DLC1) $\p(0)=0$;\\
(DLC2) $\p(a\we b)=\p(a)\we \p(b)$, for all $a,b\in A$;\\
(DLC3) If $a\in\BBBB, b\in A$ and $a\llx b$, then $(\p(a^*))^*\lle
\p(b)$;\\
(DLC4) For every $b\in\BBBB\ap$ there exists $a\in\BBBB$ such that
$b\le\p(a)$;\\
\noindent(DLC5) $\p(a)=\bigvee\{\p(b)\st b\in\BBBB, b\llx a\}$,
for every $a\in A$;

\medskip

{\noindent}let the composition $``\diamond$" of two morphisms
$\p_1:(A_1,\rho_1,\BBBB_1)\lra (A_2,\rho_2,\BBBB_2)$ and
$\p_2:(A_2,\rho_2,\BBBB_2)\lra (A_3,\rho_3,\BBBB_3)$ of\/ $\DHLC$
be defined by the formula
\begin{equation}\label{diamcon}
\p_2\diamond\p_1 = (\p_2\circ\p_1)\cuk,
\end{equation}
 where, for every
function $\psi:(A,\rho,\BBBB)\lra (B,\eta,\BBBB\ap)$ between two
objects of\/ $\DHLC$, $\psi\cuk:(A,\rho,\BBBB)\lra
(B,\eta,\BBBB\ap)$ is defined as follows:
\begin{equation}\label{cukfcon}
\psi\cuk(a)=\bigvee\{\psi(b)\st b\in \BBBB, b\llx a\},
\end{equation}
for every $a\in A$.

(We used here the same notation as in Definition \ref{dval} for
the composition between the morphisms of the category $\DHLC$ and
for the functions of the type $\psi\cuk$ because the NC-algebras can be
regarded as those LC-algebras $(A,\rho,\BBBB)$ for which $A=\BBBB$, and
hence the right sides of the formulae (\ref{cukfcon}) and
(\ref{cukfc}) coincide in the case of NC-algebras.)
\end{defi}

It can be shown that  condition (DLC3) in Definition \ref{dhc}
can be replaced by any of the following four constrains:

\medskip

\noindent(DLC3$'$) If $a,b\in\BBBB$ and $a\llx b$, then
$(\p(a^*))^*\lle \p(b)$.
\noindent(DLC3S) If $a,b\in A$ and $a\llx b$, then
$(\p(a^*))^*\lle \p(b)$.
\noindent(LC3) If, for $i=1,2$, $a_i\in\BBBB$, $b_i\in A$  and
$a_i\llx b_i$, then $\p(a_1\vee a_2)\lle \p(b_1)\vee \p(b_2)$.
\noindent(LC3S) If, for $i=1,2$, $a_i,b_i\in A$  and $a_i\llx
b_i$, then $\p(a_1\vee a_2)\lle \p(b_1)\vee \p(b_2)$.

\begin{theorem}\label{lccont}{\rm \cite{DI}}
The categories $\HLC$ and\/ $\DHLC$ are dually equi\-valent.
\end{theorem}

 The duality, constructed in Theorem \ref{lccont} and denoted by $\Psi^t:\HLC\lra\DHLC$,
 is an extension of the Roeper's correspondence $\Psi^t$ defined by (\ref{psit1})
(i.e. the definition of the contravariant functor $\Psi^t$   {\em on the objects of the category}\/ $\HLC$
coincides with the definition of the Roeper's correspondence).

We will also need a lemma from \cite{CNG}:

\begin{lm}\label{isombool}
Let $X$ be a dense subspace of a topological space $Y$. Then the
functions $r:RC(Y)\lra RC(X)$, $F\mapsto F\cap X$, and $e:RC(X)\lra
RC(Y)$, $G\mapsto \cl_Y(G)$, are Boolean isomorphisms between Boolean
algebras $RC(X)$ and $RC(Y)$, and $e\circ r=id_{RC(Y)}$, $r\circ
e=id_{RC(X)}$.
\end{lm}

For the notions and  notation not defined here see \cite{AHS, J,
E, Si}.

%--------------------------------------------------------------------
%----------------------------------------------------------------------
%%%-----------------------------------------------------------------

\section{Sums in the categories $\DHLC$ and $\DHC$}

In \cite{DII}, we described the $\DHLC$-products of complete local contact algebras. Here we will describe
the  $\DHLC$-sums of finite families of complete local contact algebras and the $\DHC$-sums of arbitrarily many complete contact algebras using the notion of a {\em sum of a family of Boolean algebras} (see \cite{GH}) which is known also as a
{\em free product} (see \cite{kop89}). (We will denote the sum of a family $\{A_\g\st\g\in\GA\}$ of Boolean algebras by $\bigoplus_{\g\in\GA} A_\g$ (as in \cite{kop89}).) Note that the sums (resp.,
finite sums) in the category $\DHC$ (resp., $\DHLC$) surely exist because the dual
category $\HC$ (resp., $\HLC$) of all compact (resp., locally compact) Hausdorff spaces and
continuous maps  has products (resp., finite products).

Let us recall the definition of the notion of a sum of a family $(A_i)_{i\in I}$ of Boolean algebras (see, e.g. \cite{kop89}): a pair $(A, (e_i)_{i\in I})$ is a {\em sum of} $(A_i)_{i\in I}$ if
$A$ is a Boolean algebra, each $e_i$ is a homomorphism from $A_i$ into $A$ and, for every family $(f_i)_{i\in I}$ of homomorphisms from $A_i$ into any Boolean algebra $B$, there is a unique homomorphism $f:A\lra B$ such that $f\circ e_i=f_i$ for $i\in I$. It is well known that every family of Boolean algebras has, up to isomorphism, a unique sum. Recall, as well, that a family $(B_i)_{i\in I}$ of subalgebras of a Boolean algebra $A$ is {\em independent}\/ if, for arbitrary $n\in\NNNN^+$, pairwise distinct $i(1),\ldots,i(n)\in I$ and non-zero elements $b_{i(k)}$ of $B_{i(k)}$, for $k=1,\ldots,n$, $b_{i(1)}\we\ldots\we b_{i(n)}>0$ in $A$. The following characterization of the sums holds (see, e.g., \cite{kop89}):

\begin{pro}\label{boosum}
Let $A$ be a Boolean algebra and, for $i\in I$, $e_i:A_i\lra A$ a homomorphism; assume that no $A_i$ is trivial. The pair $(A,(e_i)_{i\in I})$ is a sum of $(A_i)_{i\in I}$ iff each of (a) through (c) holds:

(a) each $e_i:A_i\lra A$ is an injection,

(b) $(e_i(A_i))_{i\in I}$ is an independent family of subalgebras of $A$,

(c) $A$ is generated by $\bigcup_{i\in I}e_i(A_i)$.\\
Moreover, if $(A,(e_i)_{i\in I})$ is a sum of $(A_i)_{i\in I}$ then

(d) $e_i(A_i)\cap e_j(A_j)=\{0,1\}$, for $i\neq j$.
\end{pro}

We start with a proposition which should be known, although I was not able to find it in the literature. Recall that a topological space $X$ is called {\em semiregular} if $RO(X)$
is a base of $X$. By a {\em completion} of a Boolean algebra $A$, we will understand the {\em MacNeille completion} of $A$.
% (i.e. the {\em minimal completion} of $A$).

\begin{pro}\label{rcsemi}
Let $\{X_\g\st\g\in\GA\}$ be a family of semiregular topological spaces and $X=\prod\{X_\g\st\g\in\GA\}$. Then the Boolean algebra $RC(X)$ is isomorphic to the completion of\/ $\bigoplus_{\g\in\GA} RC(X_\g)$.
\end{pro}

\doc Let, for every $\g\in\GA$, $\pi_\g:X\lra X_\g$ be the projection.  Using the fact that $\pi_\g$ is an open map (and, thus, the formulae $\cl(\pi_\g\inv(M))=\pi_\g\inv(\cl(M))$ and $\int(\pi_\g\inv(M))=\pi_\g\inv(\int(M))$ hold for every $M\sbe X_\g$)
(see, e.g., \cite{E}),  it is easy to show, that the map $\p_\g: RC(X_\g)\lra RC(X), \ F\mapsto\pi_\g\inv(F),$ is  a complete monomorphism for every $\g\in\GA$. Set $A_\g=\p_\g(RC(X_\g))$, for every $\g\in\GA$, and let $A$ be the subalgebra of $RC(X)$ generated by $\bigcup\{A_\g\st\g\in\GA\}$. It is easy to check that, for every finite non-empty subset $\GA_0$ of $\GA$, we have that if $a_\g\in A_\g\stm\{0\}$ for every $\g\in\GA_0$, then $\bw\{a_\g\st\g\in\GA_0\}\neq 0$ (i.e. the family $\{A_\g\st\g\in\GA\}$ is an {\em independent family} (see, e.g., \cite{kop89})). Thus, by
%\cite[Lemma 1, page 428]{GH} or
\cite[Proposition 11.4]{kop89}, we get that $A=\bigoplus_{\g\in\GA} RC(X_\g)$. Since $RO(X_\g)$ is a base of $X_\g$, for every $\g\in\GA$, we obtain that $A$ is a dense subalgebra of $RC(X)$. Thus, $RC(X)$ is the completion of $A$.
%Finally, note that $A$ is isomorphic to $\bigoplus_{\g\in\GA} RC(X_\g)$.
\sqs

The proof of this proposition shows that the following is even true:

\begin{cor}\label{rcsemic}
Let $\{X_\g\st\g\in\GA\}$ be a family of semiregular topological spaces and $X=\prod\{X_\g\st\g\in\GA\}$. Let, for every $\g\in\GA$, $B_\g$ be a subalgebra of $RC(X_\g)$ such that $\{\int(F)\st F\in B_\g\}$ is a base of $X_\g$. Then the Boolean algebra $RC(X)$ is isomorphic to the completion of\/ $\bigoplus_{\g\in\GA} B_\g$.
\end{cor}

\begin{defi}\label{deflcasum}
\rm
Let $n\in\mathbb{N}^+$ and let, for every $\ion$, $(A_i,\rho_i,\BBBB_i)$ be a CLCA. Let $$(A,(\p_i)_{i=1}^n)=\bigoplus_{i=1}^n A_i,$$ where, for every $i\in\oon$, $$\p_i:A_i\lra A$$ is the canonical complete monomorphism, and let $\ti{A}$ be the completion of $A$. We can  suppose, without loose of generality, that $A\sbe \ti{A}$. Set $$E=\{\bw_{i=1}^n \p_i(a_i)\st a_i\in \BBBB_i\}$$ and let\/ $\wt{\BBBB}$ be the ideal of $\ti{A}$
generated by $E$ (thus, $$\wt{\BBBB}=\{x\in\ti{A}\st x\le e_1\vee \ldots\vee e_n\mbox{ for some }n\in\NNNN^+\mbox{ and }e_1,\ldots, e_n\in E\}).$$ For every two elements $a=\bw_{i=1}^n\p_i(a_i)$ and $b=\bw_{i=1}^n\p_i(b_i)$ of $E$, set $$a\ti{\rho} b\Leftrightarrow (a_i\rho_i b_i,\fa i\in\oon).$$
Further, for every two elements $c$ and $d$ of\/ $\wt{\BBBB}$, set $$c(-\ti{\rho})d\Leftrightarrow (\ex k,l\in\mathbb{N}^+ \mbox{ and }\ex c_1,\ldots,c_k,d_1,\ldots, d_l\in E
\mbox{ such that }$$ $$c\le\bv_{i=1}^k c_i,\ d\le\bv_{j=1}^l d_j\mbox{ and } c_i(-\ti{\rho})d_j,\ \fa i=1,\ldots, k\mbox{ and }\fa j=1,\ldots, l).$$ Finally, for every two elements $a$ and $b$ of $\ti{A}$, set $$a\ti{\rho}b\Leftrightarrow (\ex c,d\in\wt{\BBBB}\mbox{ such that } c\le a,\ d\le b\mbox{ and } c\ti{\rho} d).$$ Then the triple $(\ti{A},\ti{\rho},\wt{\BBBB})$ will be denoted by $\bigoplus_{i=1}^n (A_i,\rho_i,\BBBB_i)$.
\end{defi}

\begin{theorem}\label{lcasum}
Let $n\in{\NNNN^+}$ and $\{(A_i,\rho_i,\BBBB_i)\st \ion\}$ be a family of CLCAs. Then $\bigoplus_{i=1}^n (A_i,\rho_i,\BBBB_i)$ is a $\DHLC$-sum of the family $\{(A_i,\rho_i,\BBBB_i)\st \ion\}$.
\end{theorem}

\doc As the Duality Theorem \ref{lccont} shows, for every $i\in\oon$ there exists a $X_i\in|\HLC|$ such that the CLCAs $(RC(X_i),\rho_{X_i},CR(X_i))$ and
$(A_i,\rho_i,\BBBB_i)$ are LCA-isomorphic. Let $X=\prod_{i=1}^n X_i$. Then we have, in the notation of Definition \ref{deflcasum}, that the Boolean algebras $RC(X)$
and $\ti{A}$ are isomorphic (see Proposition \ref{rcsemi}). Also, again in the notation of Definition \ref{deflcasum}, $(A,(\p_i)_{i=1}^n)$ is isomorphic to
$(\bigoplus_{i=1}^n RC(X_i), (\psi_i)_{i=1}^n)$, where $\psi_i: RC(X_i)\lra RC(X),\ F\mapsto\pi_i\inv(F)$, and $\pi_i:X\lra X_i$ is the projection, for every $i\in\oon$
(this follows from Proposition \ref{boosum}).
Thus, the set $E$ from Definition \ref{deflcasum} corresponds to the following set:
 $$E'=\{\bw_{i=1}^n\psi_i(F_i)\st F_i\in CR(X_i)\}.$$ 
 Let $F\in E'$. Then there exist $F_i\in CR(X_i)$, for $\ion$, such that $F=\bw_{i=1}^n\psi_i(F_i)$. Set $U_i=\int_{X_i}(F_i)$, for $\ion$. Then $F=\bw_{i=1}^n\pi_i\inv(F_i)=\cl_X(\bigcap_{i=1}^n\int_X(\pi_i\inv(F_i)))=\cl_X(\bigcap_{i=1}^n\pi_i\inv(U_i))=\cl(\prod_{i=1}^n U_i)=\prod_{i=1}^n F_i$ (note that we used \cite[1.4.C,2.3.3]{E} here). Hence, for every $F,G\in E'$, where $F=\prod_{i=1}^n F_i$ and $G=\prod_{i=1}^n G_i$, we have that $$F\rho_X G\Leftrightarrow F\cap G\nes\Leftrightarrow ( F_i\cap G_i\nes, \fa i=1,\ldots, n)\Leftrightarrow ( F_i\rho_{X_i} G_i, \fa i=1,\ldots, n).$$ Further, since $\{\prod_{i=1}^n U_i \st  U_i\in RO(X_i), \fa \ion\}$ is a base of $X$ and $X$ is regular, we obtain that $CR(X)$ coincides with the ideal of $RC(X)$ generated by $E'$. The fact that every two disjoint compact subsets of $X$ can be separated by open sets implies that if $F,G\in CR(X)$ then $F(-\rho_X) G$ (i.e. $F\cap G=\ems$) iff there exists finitely many elements $F_1,\ldots,F_k,G_1,\ldots,G_l\in E'$ such that $F\sbe\bigcup_{i=1}^k F_i$,  $G\sbe\bigcup_{i=1}^l G_i$ and $F_i\cap G_j=\ems$ (i.e. $F_i(-\rho_X) G_j$) for all $i=1,\ldots,k$ and all $j=1,\ldots,l$. Finally, since $(RC(X),\rho_X,CR(X))$ is an LCA (see \ref{stanlocn}), we have (by (BC2)) that for any $F',G'\in RC(X)$, $F'\rho_X G'\Leftrightarrow \ex F,G\in CR(X)$ such that $F\sbe F'$, $G\sbe G'$ and $F\rho_X G$. All this shows that  the triple $(\ti{A},\ti{\rho},\wt{\BBBB})$ from \ref{deflcasum} is an LCA which is LCA-isomorphic to $(RC(X),\rho_X,CR(X))$. Now, using Theorem \ref{lccont} and the facts that
 $\Psi^t(X)=(RC(X),\rho_X,CR(X))$, $\Psi^t(X_i)=(RC(X_i),\rho_{X_i},CR(X_i))$ for all $\ion$, and $X$ is a $\HLC$-product of the family $\{X_i\st\ion\}$, we get that
 $(RC(X),\rho_X,CR(X))$ is a $\DHLC$-sum of the family $\{(RC(X_i),\rho_{X_i},CR(X_i))\st \ion\}$. Thus $(\ti{A},\ti{\rho},\wt{\BBBB})$ is a $\DHLC$-sum of the family $\{(A_i,\rho_i,\BBBB_i)\st \ion\}$. \sqs

\begin{defi}\label{defncasum}
\rm
Let $J$ be a set and let, for every $j\in J$, $(A_j,\rho_j)$ be a CNCA. Let $$(A,(\p_j)_{j\in J})=\bigoplus_{j\in J} A_j,$$ where, for every $j\in J$, $$\p_j:A_j\lra A$$ is the canonical complete monomorphism, and let $\ti{A}$ be the completion of $A$. We can  suppose, without loose of generality, that $A\sbe \ti{A}$. Set
$$E=\{\bw_{i\in I}\p_i(a_i)\st  I\sbe J, |I|<\aleph_0,  a_i\in A_i, \fa i\in I \}.$$  For every two elements $a=\bw_{i\in I_1}\p_i(a_i)$ and $b=\bw_{i\in I_2}\p_i(b_i)$ of $E$, set $$a\ti{\rho} b\Leftrightarrow ( a_i\rho_i b_i,\fa i\in I_1\cap I_2).$$
Further, for every two elements $c$ and $d$ of $\ti{A}$, set $$c(-\ti{\rho})d\Leftrightarrow (\ex k,l\in\mathbb{N}^+ \mbox{ and }\ex c_1,\ldots,c_k,d_1,\ldots, d_l\in E
\mbox{ such that }$$ $$c\le\bv_{i=1}^k c_i,\ d\le\bv_{j=1}^l d_j\mbox{ and } c_i(-\ti{\rho})d_j,\ \fa i=1,\ldots, k\mbox{ and }\fa j=1,\ldots, l).$$  Then the pair $(\ti{A},\ti{\rho})$ will be denoted by $\bigoplus_{j\in J} (A_j,\rho_j)$.
\end{defi}

\begin{theorem}\label{ncasum}
Let  $\{(A_j,\rho_j)\st j\in J\}$ be a family of CNCAs. Then $\bigoplus_{j\in J} (A_j,\rho_j)$ is a $\DHC$-sum of the family $\{(A_j,\rho_j)\st j\in J\}$.
\end{theorem}

\doc The proof is similar  to that one of Theorem \ref{lcasum}. In it de Vries' Duality Theorem \ref{dvth} instead of Theorem \ref{lccont} can be used. \sqs

\section{A Whiteheadian-type description of Euclidean spaces}

\begin{nota}\label{notation}
\rm
We will denote by $\mathbb{Z}$ the set of all integers with the natural order,  by $\mathbb{I}$ the unit interval $[0,1]$ with its natural topology
and by $\mathbb{I}\ap$  -- the open interval $(0,1)$ with its natural topology, by $\NNNN$ the set of natural numbers, by $\mathbb{J}$ the subspace of the real line
consisting of all irrational numbers,  and by $\mathbb{D}$ the set of all dyadic numbers in the
 interval $(0,1)$. We set $\mathbb{Z}_0=\mathbb{Z}\stm\{0\}$, $\mathbb{Z}^-=\mathbb{Z}\stm\NNNN$ and $\mathbb{J}_2=\mathbb{I}\ap\stm\mathbb{D}$.
 If $(X,<)$ is a linearly ordered set and $x\in X$, then we set
$$ succ(x)=\{y\in X\st x< y\},\  pred(x)=\{y\in X\st y< x\};$$ also, we denote by $x^+$  the {\em successor of}\/ $x$ (when it exists) and by $x^-$ -- the {\em predecessor
of}\/ $x$ (when it exists).
If $M$ is a set, then we will denote by $P(M)$ the power set Boolean algebra  of $M$;  the cardinality of $M$ will be denoted by $|M|$. If $X$ is a topological
space, then we will denote by $CO(X)$ the set of all clopen (= closed and open) subsets of $X$.
\end{nota}

Now we will construct a CLCA $(\ti{A},\ts,\wt{\BBBB})$ and we will show that it is LCA-isomorphic to $\Psi^t(\mathbb{R})$.

\begin{nist}\label{constr}{\bf The construction of $(\ti{A},\ts,\wt{\BBBB})$.}
\rm
Let $A_i=P(\mathbb{Z}_0)$, for every $i\in{\NNNN^+}$. Thus, if $i\in{\NNNN^+}$ and $a_i\in A_i$, then $a_i$ is a subset of $\mathbb{Z}_0$ and its cardinality will be denoted by $|a_i|$.  Let $(A,(\p_i)_{i\in{\NNNN^+}})$ be the sum of
 Boolean algebras $\{A_i\st i\in{\NNNN^+}\}$; then, by Proposition \ref{boosum}, for every $i\in{\NNNN^+}$, $\p_i:A_i\lra A$ is a monomorphism, 
 the family $\{\p_i(A_i)\st i\in{\NNNN^+}\}$ is an  independent family  and the set $\bigcup_{i\in{\NNNN^+}}\p_i(A_i)$ generates $A$.
%$A=\bigoplus_{i\in{\NNNN^+}} A_i$ and
Let $\ti{A}$ be the completion of $A$.
We can  suppose, without loose of generality, that $A\sbe \ti{A}$.
%and $A_i\sbe A$, for all $i\in{\NNNN^+}$.
%If $i\in{\NNNN^+}$ and $a\in P(\mathbb{Z}_0)$, we will sometimes write $``(a)_i$" when we want to stress on the fact that  $a$ is regarded as an element of $A_i$.

The following subset of $A$ will be important for us:
\begin{equation}\label{B0}
B_0=\{\p_1(a_1)\we \ldots\we \p_k(a_k)\st k\in{\NNNN^+}, (\fa i=1,\ldots, k)(a_i\in A_i \mbox{ and } |a_i|=1)\}.
\end{equation}
If $b\in B_0$ and $b=\p_1(a_1)\we \ldots\we \p_k(a_k)$, where $a_k=\{p\}$, then we set
\begin{equation}\label{bu}
b_-=\p_1(a_1)\we \p_2(a_2)\we \ldots\we \p_{k-1}(a_{k-1})\we \p_k(\{p^-\}).
\end{equation}
For every $b\in B_0$, where $b=\p_1(a_1)\we \ldots\we \p_k(a_k)$, and every $n\in{\NNNN^+}$, we set
\begin{equation}\label{qbn}
q_{bn}=(b_-\we\p_{k+1}( succ(n)))\vee(b\we\p_{k+1}(pred(-n))).
\end{equation}
Now we set
\begin{equation}\label{B1}
B_1=\{q_{bn}\st b\in B_0,n\in{\NNNN^+}\}.
\end{equation}

Let $\wt{\BBBB}$ be the ideal of $\ti{A}$ generated by the set $B_0\cup B_1$.
Now, we will define a relation $\ti{\s}$ on $\ti{A}$. It will be, by definition, a symmetric relation.

Let  $r,r\ap\in{\NNNN^+}$, $b,b\ap\in B_0$, $b=\p_1(a_1)\we \ldots\we \p_k(a_k)$,  $b\ap=\p_1(a_1\ap)\we \ldots\we \p_l(a_l\ap)$ and  $a_k=\{n\}$,  $a_k\ap=\{m\}$. We can suppose,
 without loose of generality, that $k\le l$. If $k<l$, then let $a_{k+1}\ap=\{p\}$.  Now we set
%$$
\begin{equation}\label{bbap}
b\ts b\ap\Leftrightarrow
%\mbox{ iff }
[(a_i=a_i\ap,\ \fa i\in\{1,\ldots, k-1\})
\end{equation}
$$\mbox{ and }(\left\{
\begin{array}{ll}
m\in\{n^-,n,n^+\}, & \mbox{if $k=l$}\\
m=n, & \mbox{if $k< l$}
\end{array}
)],\right.
$$
and
 %$$
 \begin{equation}\label{ububap}
 q_{br}\ts q_{b\ap r\ap}\Leftrightarrow [(a_i=a_i\ap,\ \fa i\in\{1,\ldots, k-1\})\mbox{ and }
 \end{equation}
 $$( \left\{
 \begin{array}{lll}
 m=n,
 %a_i=a_i\ap,\ \fa i\in\{1,\ldots, k\},
 & \mbox{if $l=k$}\\
 %(a_i=a_i\ap,\
 %\fa i\in\{1,\ldots, k-1\})
 %\mbox{ for }
 %i=1\div k-1)
 %\mbox{ and }
 %\ \& \
 (m=n
 %\ \& \
 \mbox{ and } p\le -r)\mbox{ or }(m=n^-
 %\ \& \
 \mbox{ and } p> r), & \mbox{if $l=k+1$}\ )].\\
 %(a_i=a_i\ap,\ \fa i\in\{1,\ldots, k\})\mbox{ and } (p < -r), & \mbox{if $l > k+1.$}
 (m=n
 %\ \& \
 \mbox{ and } p< -r)\mbox{ or }(m=n^-
 %\ \& \
 \mbox{ and } p > r), & \mbox{if $l > k+1$}
\end{array}\right.
%\end{equation}
$$

Let  $r\in{\NNNN^+}$, $b,b\ap\in B_0$, $b=\p_1(a_1)\we \ldots\we \p_k(a_k)$,  $b\ap=\p_1(a_1\ap)\we \ldots\we \p_l(a_l\ap)$ and  $a_k=\{n\}$,  $a_k\ap=\{m\}$.
If $k<l$, then let $a_{k+1}\ap=\{p\}$. Now, if $k> l$,  we set
%$$
\begin{equation}\label{ubbap}
q_{br}\ts b\ap\Leftrightarrow  (a_i=a_i\ap,\ \fa i\in\{1,\ldots, l\});
\end{equation}
%$$
if $k\le l$,  we set
\begin{equation}\label{ubrbap}
q_{br}\ts b\ap  \Leftrightarrow  [(a_i=a_i\ap,\ \fa i\in\{1,\ldots, k-1\})\mbox{ and }
\end{equation}
$$(
    \left\{
 \begin{array}{lll}
 m\in\{n^-,n\}, & \mbox{if $l=k$}\\
 (p\ge r\mbox{ and $m=n^-$) or } (p\le -r \mbox{ and $m=n$}), & \mbox{if $l=k+1$}\\
 (p> r\mbox{ and $m=n^-$) or } (p< -r \mbox{ and $m=n$}), & \mbox{if $l > k+1$}
\end{array}
)].\right.
$$
Further, for every two elements $c$ and $d$ of\/ $\wt{\BBBB}$, set
%$$
\begin{equation}\label{csd}
c(-\ts)d\Leftrightarrow (\ex k,l\in\mathbb{N}^+ \mbox{ and }\ex c_1,\ldots,c_k,d_1,\ldots, d_l\in B_0\cup B_1
\mbox{ such that }
\end{equation}
%$$
$$c\le\bv_{i=1}^k c_i,\ d\le\bv_{j=1}^l d_j\mbox{ and } c_i(-\ts)d_j,\ \fa i=1,\ldots, k\mbox{ and }\fa j=1,\ldots, l).$$
Finally, for every two elements $a$ and $b$ of $\ti{A}$, set
%$$
\begin{equation}\label{asb}
a\ts b\Leftrightarrow (\ex c,d\in\wt{\BBBB}\mbox{ such that } c\le a,\ d\le b\mbox{ and } c\ts d).
\end{equation}
%$$
\end{nist}

\begin{theorem}\label{whr}
The triple $(\ti{A},\ts,\wt{\BBBB})$, constructed in \ref{constr}, is a CLCA; it is LCA-isomorphic to the CLCA $(RC(\mathbb{R}),\rho_\mathbb{R},CR(\mathbb{R}))$.
Thus, the triple $(\ti{A},\ts,\wt{\BBBB})$ completely determines the real line $\mathbb{R}$ with its natural topology.
\end{theorem}

\doc In this proof, we will use the notation introduced in \ref{constr}.

Let $\mathbb{Z}_0$ be endowed with the discrete topology. Then $RC(\mathbb{Z}_0)=P(\mathbb{Z}_0)$ and Proposition \ref{rcsemi} shows that the algebra $\ti{A}$, constructed in \ref{constr}, is isomorphic to $RC(\mathbb{Z}_0^{{\NNNN^+}})$. Since the space $\mathbb{Z}_0^{{\NNNN^+}}$ is homeomorphic to $\mathbb{J}$ (see, e.g., \cite{E}), we get, by Lemma \ref{isombool}, that $\ti{A}$ is isomorphic to $RC(\mathbb{R})$. Clearly, $RC(\mathbb{J})$ can be endowed with an LCA-structure LCA-isomorphic to the LCA $(RC(\mathbb{R}),\rho_\mathbb{R},CR(\mathbb{R}))$. Then, using the homeomorphism between  $\mathbb{J}$  and $\mathbb{Z}_0^{{\NNNN^+}}$, we can transfer this structure to $RC(\mathbb{Z}_0^{{\NNNN^+}})$ and, hence, to $\ti{A}$. For technical reasons, this plan will be slightly modified. We will use the homeomorphism between
$\mathbb{Z}_0^{{\NNNN^+}}$ and $\mathbb{J}_2$ described in \cite{Alex}. Since $\mathbb{J}_2$ is dense in the open interval $\mathbb{I}\ap$, and $\mathbb{I}\ap$ is homeomorphic to $\mathbb{R}$, we can use $\mathbb{J}_2$ instead of $\mathbb{J}$ for realizing the desired transfer. So, we start with the description (given by P. S. Alexandroff \cite{Alex}) of the homeomorphism $f:\mathbb{Z}_0^{{\NNNN^+}}\lra\mathbb{J}_2$. Let, for every $j\in\mathbb{N}^+$, $\DE_j=[1-\frac{1}{2^j},1-\frac{1}{2^{j+1}}]$ and let, for every $j\in\mathbb{Z}^-$, $\DE_j=[2^{j-1},2^j]$. Set $\d_1=\{\DE_j\st j\in\mathbb{Z}_0\}$. Further, for every $\DE_j\in\d_1$, where $\DE_j=[a_j,b_j]$, set $d_j=b_j-a_j$ and $\DE_{jk}=[b_j-\frac{d_j}{2^k},b_j-\frac{d_j}{2^{k+1}}]$ when $k\in\mathbb{N}^+$,
$\DE_{jk}=[a_j+d_j.2^{k-1},a_j+d_j.2^k]$ when $k\in\mathbb{Z}^-$. Let $\d_2=\{\DE_{jk}\st j,k\in\mathbb{Z}_0\}$.  In the next step we construct analogously the family $\d_3$, and so on.
Set $\d=\bigcup\{\d_i\st i\in{\NNNN^+}\}$.
It is easy to see that the set of all end-points of the elements of the family $\d$ coincides with the set $\mathbb{D}$. Now we define the function $f:\mathbb{Z}_0^{{\NNNN^+}}\lra \mathbb{J}_2$ by the formula $$f(n_1,n_2,\ldots,n_k,\ldots)=\DE_{n_1}\cap\DE_{n_1 n_2}\cap\ldots\cap\DE_{n_1 n_2\ldots n_k}\cap\ldots .$$
One can prove that the definition of $f$ is correct and that $f$ is a homeomorphism. Set $X_i=\mathbb{Z}_0$, for every $i\in{\NNNN^+}$. Let $X=\prod\{X_i\st i\in{\NNNN^+}\}$ and let $$\pi_i:X\lra X_i,$$ where $i\in{\NNNN^+}$, be the projection. Then, for every $k\in{\NNNN^+}$ and every $n_i\in X_i$, where $i=1,\ldots,k$, we have that (writing, for short, $``\pi_i\inv(n_i)$" instead of $``\pi_i\inv(\{n_i\})$")
\begin{equation}\label{deltak}
f(\bigcap_{i=1}^k\pi_i\inv(n_i))=\DE_{n_1 n_2\ldots n_k}\cap\mathbb{J}_2.
\end{equation}

Let $\psi_i:RC(X_i)\lra RC(X)$, $F\mapsto\pi_i\inv(F)$, where $i\in{\NNNN^+}$; then, as we have seen in the proof of Proposition \ref{rcsemi}, $\psi_i$ is a complete monomorphism. Set $A_i\ap=\psi_i(RC(X_i))$. Since $X_i$ is a discrete space, we have that $A_i=RC(X_i)$ and $A_i\ap\sbe CO(X)$,  for all $i\in{\NNNN^+}$. Thus, for the elements of the subset $\bigcup_{i\in{\NNNN^+}}A_i\ap$ of $RC(X)$, the Boolean operation $``$meet in $RC(X)$"
coincides with the set-theoretic operation $``$intersection" between the subsets of $X$, and the same for the Boolean complement in $RC(X)$ and the set-theoretic complement in $X$. We also have that the Boolean algebras $A_i$ and $A_i\ap$ are isomorphic. Let $A\ap$ be the subalgebra of $P(X)$ generated by $\bigcup_{i\in{\NNNN^+}}A_i\ap$. Then $A\ap$ is isomorphic to $A$. Note that $A\ap$ is a subalgebra of $CO(X)$. Also, $A\ap$ is a dense subalgebra of $RC(X)$; therefore, $RC(X)$ is the completion  of $A\ap$. Thus, $\ti{A}$ is isomorphic to $RC(X)$. So, without loose of generality, we can think that $\ti{A}$ is $RC(X)$, $A$ is $A\ap$, $\p_i=\psi_i$ and hence $\p_i(A_i)$ is $A_i\ap$, for $i\in{\NNNN^+}$.
  We will now construct an  LCA $(RC(X),\s,\BBBB)$   LCA-isomorphic to $(RC(\mathbb{R}),\rho_{\mathbb{R}},CR(\mathbb{R}))$. Then, identifying $RC(X)$ with $\ti{A}$, we will
  show that $\s=\ts$ and $\BBBB=\wt{\BBBB}$.

Let $\BBBB_2=\{M\in RC(\mathbb{J}_2)\st \cl_{\mathbb{I}\ap}(M)$ is compact$\}$. For every two elements $M$ and $N$ of $RC(\mathbb{J}_2)$, set $M\rho_2 N\Leftrightarrow
\cl_{\mathbb{I}\ap}(M)\cap\cl_{\mathbb{I}\ap}(N)\nes$. Then, using Lemma \ref{isombool}, we get that the triple $(RC(\mathbb{J}_2),\rho_2,\BBBB_2)$ is LCA-isomorphic to the LCA $(RC(\mathbb{I}\ap),\rho_{\mathbb{I}\ap},CR(\mathbb{I}\ap))$ (which, in turn, is LCA-isomorphic to $(RC(\mathbb{R}),\rho_{\mathbb{R}},CR(\mathbb{R}))$).
Now, for every two elements $F,G\in RC(X)$, we set
%$$
\begin{equation}\label{sigrho}
F\s G\Leftrightarrow f(F)\rho_2 f(G).
\end{equation}
%$$
Also, we put
\begin{equation}\label{BBBB}
\BBBB=\{f\inv(M)\st M\in\BBBB_2\}.
 \end{equation}
 Obviously, $(RC(X),\s,\BBBB)$  is LCA-isomorphic to $(RC(\mathbb{R}),\rho_{\mathbb{R}},CR(\mathbb{R}))$. In the rest of this proof, we will show that the definitions of $\BBBB$ and $\s$ given above agree with the corresponding definitions of $\wt{\BBBB}$ and $\ti{\s}$ given in \ref{constr}.

Note first that the subset $B_0\ap$ of $A\ap$, which corresponds to the subset $B_0$ of $A$ described in \ref{constr}, is the following:
\begin{equation}\label{B0ap}
B_0\ap=\{ \bigcap_{i=1}^k\pi_i\inv(n_i)\st k\in{\NNNN^+},(\fa i=1,\ldots, k)(n_i\in X_i)\}.
\end{equation}
Let $F,G\in B_0\ap$ and $F=\bigcap_{i=1}^k\pi_i\inv(n_i)$, $G=\bigcap_{i=1}^l\pi_i\inv(m_i)$. We can suppose,
 without loose of generality, that $k\le l$.
% and set $n=n_k$, $m=m_k$.
Then, by (\ref{deltak}) and Lemma \ref{isombool},
$\cl_{\mathbb{I}\ap}(f(F))=\DE_{n_1 n_2\ldots n_k}$ and  $\cl_{\mathbb{I}\ap}(f(G))=\DE_{m_1 m_2\ldots m_l}$. If $k=l$, then, clearly,
 $\DE_{n_1 n_2\ldots n_k}\cap\DE_{m_1 m_2\ldots m_k}\nes$ iff ($n_i=m_i$, for all $i=1,\ldots, k-1$, and $m_k\in\{n_k^-,n_k,n_k^+\}$).
 If $k< l$, then, obviously, $\DE_{n_1 n_2\ldots n_k}\cap\DE_{m_1 m_2\ldots m_l}\nes$ iff ($n_i=m_i$, for all $i=1,\ldots, k$).
 %(i.e. iff ($n_i=m_i$, for all $i=1,\ldots, k-1$, and $m=n$)).
 Then, using (\ref{sigrho}) and the formula (\ref{bbap}), we get that  $\s$ and $\ts$ agree on $B_0\ap$ (or, equivalently, on $B_0$).

Let $F\in B_0\ap$,  $F=\bigcap_{i=1}^k\pi_i\inv(n_i)$ and $n\in{\NNNN^+}$. Then the element $Q_{Fn}$ of $A\ap$ corresponding to the element $q_{bn}$ of $A$, where $b\in B_0$ corresponds to $F$,
is the following: $$Q_{Fn}=[(\bigcap_{i=1}^{k-1}\pi_i\inv(n_i))\cap\pi_k\inv(n_k^-)\cap\pi_{k+1}\inv( succ(n))]\cup[F\cap\pi_{k+1}\inv(pred(-n))].$$
 Clearly,
 \begin{equation}\label{QFn}
 Q_{Fn}=
 [\bigcup_{s\in succ(n)}(\bigcap_{i=1}^{k-1}\pi_i\inv(n_i)\cap\pi_k\inv(n_k^-)\cap\pi_{k+1}\inv(s))]\ \cup
 \end{equation}
 $$[\bigcup_{s\in pred(-n)}(\bigcap_{i=1}^{k}\pi_i\inv(n_i)\cap\pi_{k+1}\inv(s))].$$
 (It is easy to see, as well, that in the formula (\ref{QFn}) the sign of the union can be replaced everywhere with the sign of the join in $RC(X)$.)
 Thus,
 \begin{equation}\label{UFn}
 f(Q_{Fn})=[(\bigcup_{s\in succ(n)}\DE_{n_1 n_2\ldots n_{k-1}n_k^-s})\ \cup\ (\bigcup_{s\in pred(-n)}\DE_{n_1 n_2\ldots n_{k}s})]\cap\mathbb{J}_2.
 \end{equation}
 Let $d$ be the left end-point of the closed interval $\DE_{n_1 n_2\ldots n_{k}}$. Then it is easy to see that
 \begin{equation}\label{UFnc}
 \cl_{\mathbb{I}\ap}(f(Q_{Fn}))=[d-\ep_n,d+\ep_n\ap],
  \end{equation}
  where $\ep_n$ and $\ep_n\ap$ depend  from $n$ and also from $n_1,\ldots,n_k$ (for simplicity, we don't reflect this dependence on the notation), but for fixed  $n_1,\ldots,n_k$, we have that $\ep_n>\ep_{n+1}>0$, $\ep_n\ap>\ep_{n+1}\ap>0$, for all $n\in{\NNNN^+}$, and
  $\lim_{n\rightarrow\infty}\ep_n=0$, $\lim_{n\rightarrow\infty}\ep_n\ap=0$; also, the closed interval $[d-\ep_n,d+\ep_n\ap]$ lies in the open interval having as end-points the middles of the closed intervals $\DE_{n_1 n_2\ldots n_{k-1}n_k^-}$ and $\DE_{n_1 n_2\ldots n_{k}}$.
  Since the family $\{D\cap\mathbb{J}_2\st D\in\d\}$ is a base of $\mathbb{J}_2$ and every element of $\mathbb{D}$ appears as a left end-point of some element of the family $\d$, we get that the family
  $$\BB=\{\int_{\mathbb{I}\ap}(\cl_{\mathbb{I}\ap}((f(F))),\int_{\mathbb{I}\ap}(\cl_{\mathbb{I}\ap}((f(Q_{Fn})))\st n\in{\NNNN^+}, F\in B_0\ap\}$$ is a base of $\mathbb{I}\ap$.
Also, if $$\B=\{\cl_{\mathbb{I}\ap}((f(F)),\cl_{\mathbb{I}\ap}((f(Q_{Fn}))\st n\in{\NNNN^+}, F\in B_0\ap\},$$ then $\B=\{\cl_{\mathbb{I}\ap}(U)\st U\in\BB\}$ and $\B\sbe CR(\mathbb{I}\ap)$. Hence, $\B$ generates the ideal $CR(\mathbb{I}\ap)$ of $RC(\mathbb{I}\ap)$.
Clearly, the family
\begin{equation}\label{B1ap}
B_1\ap=\{Q_{Fn}\st F\in B_0\ap, n\in{\NNNN^+}\}
\end{equation}
 corresponds to the subset $B_1$ of $A$ constructed in \ref{constr}. Since $\B=\{\cl_{\mathbb{I}\ap}(G)\st G\in f(B_0\ap\cup B_1\ap)\}$, we get that the subset $f(B_0\ap\cup B_1\ap)$ of $RC(\mathbb{J}_2)$ generates the ideal $\BBBB_2$ of $RC(\mathbb{J}_2)$. Thus, the subset
$B_0\ap\cup B_1\ap$ of $RC(X)$ generates the ideal $\BBBB$ of $RC(X)$. Therefore, $\BBBB$ corresponds to $\wt{\BBBB}$; we can even write that $\BBBB=\wt{\BBBB}$.

Let now   $r,r\ap\in{\NNNN^+}$, $F,F\ap\in B_0\ap$, $F=\pi_1\inv(n_1)\cap \ldots\cap \pi_k\inv(n_k)$ and  $F\ap=\pi_1\inv(n_1\ap)\cap \ldots\cap \pi_l\inv(n_l\ap)$.
 We can suppose,
 without loose of generality, that $k\le l$.
 Let $d$ and $d\ap$ be the left end-points of the closed intervals $\DE_{n_1 n_2\ldots n_{k}}$ and $\DE_{n_1\ap n_2\ap\ldots n_{l}\ap}$, respectively.
 Then, using (\ref{UFnc}), we get that $\cl_{\mathbb{I}\ap}(f(Q_{Fr}))=[d-\ep_r,d+\ep_r\ap]$ and $\cl_{\mathbb{I}\ap}(f(Q_{F\ap r\ap}))=[d\ap-\ep_{r\ap},d\ap+\ep_{r\ap}\ap]$.
 If $k=l$, then it is easy to see that $\cl_{\mathbb{I}\ap}(f(Q_{Fr}))\cap\cl_{\mathbb{I}\ap}(f(Q_{F\ap r\ap}))\nes$ iff ($n_i=n_i\ap$, for all $i=1,\ldots,k$). If $l=k+1$, then one readily checks that $\cl_{\mathbb{I}\ap}(f(Q_{Fr}))\cap\cl_{\mathbb{I}\ap}(f(Q_{F\ap r\ap}))\nes$ iff [($n_i=n_i\ap$, for all $i=1,\ldots,k-1$) and (($n_k=n_k\ap$ and $n_{k+1}\ap\le -r$) or ($n_k\ap=(n_k)^-$ and $n_{k+1}\ap > r$))]. Finally, if $l>k+1$, then $\cl_{\mathbb{I}\ap}(f(Q_{Fr}))\cap\cl_{\mathbb{I}\ap}(f(Q_{F\ap r\ap}))\nes$ iff [($n_i=n_i\ap$, for all $i=1,\ldots,k-1$) and
 (($n_k=n_k\ap$ and $n_{k+1}\ap < -r$) or ($n_k\ap=(n_k)^-$ and $n_{k+1}\ap > r$))].
  All this shows that the relations $\s$ and $\ts$ agree on $B_1\ap$ (or, equivalently, on $B_1$).

 Let    $r\in{\NNNN^+}$, $F,F\ap\in B_0\ap$, $F=\pi_1\inv(n_1)\cap \ldots\cap \pi_k\inv(n_k)$ and  $F\ap=\pi_1\inv(n_1\ap)\cap \ldots\cap \pi_l\inv(n_l\ap)$.
 If $l<k$, then we get that $\cl_{\mathbb{I}\ap}(f(Q_{Fr}))\cap\cl_{\mathbb{I}\ap}(f(F\ap))\nes$ iff ($n_i=n_i\ap$, for all $i=1,\ldots,l$). If $l=k$, then
 $\cl_{\mathbb{I}\ap}(f(Q_{Fr}))\cap\cl_{\mathbb{I}\ap}(f(F\ap))\nes$ iff ($n_i=n_i\ap$, for all $i=1,\ldots,k-1$, and $n_k\ap\in\{n_k^-,n_k\}$). If $l=k+1$, then
 $\cl_{\mathbb{I}\ap}(f(Q_{Fr}))\cap\cl_{\mathbb{I}\ap}(f(F\ap))\nes$ iff [($n_i=n_i\ap$, for all $i=1,\ldots,k-1$), and (($n_k\ap=n_k^-$ and $n_{k+1}\ap\ge r$) or ($n_k\ap=n_k$ and $n_{k+1}\ap\le -r$))]. Finally, if $l> k+1$, then
 $\cl_{\mathbb{I}\ap}(f(Q_{Fr}))\cap\cl_{\mathbb{I}\ap}(f(F\ap))\nes$ iff [($n_i=n_i\ap$, for all $i=1,\ldots,k-1$), and (($n_k\ap=n_k^-$ and $n_{k+1}\ap > r$) or ($n_k\ap=n_k$ and $n_{k+1}\ap < -r$))]. We get that the relations $\s$ and $\ts$ agree on $B_0\ap\cup B_1\ap$ (or, equivalently, on $B_0\cup B_1$).

 Now, using the facts that $\BB$ is a base of $\mathbb{I}\ap$,  $\mathbb{I}\ap$ is a regular space, and $\cl_{\mathbb{I}\ap}(f(F))$ is a compact set for all $F\in\BBBB$,
 we get that for all $F,G\in\BBBB$, $\cl_{\mathbb{I}\ap}(f(F))\cap \cl_{\mathbb{I}\ap}(f(G))=\ems$ iff (there exist $F_1,\ldots,F_k,G_1,\ldots,G_l\in B_0\ap\cup B_1\ap$ such that $F\sbe\bigcup_{i=1}^k F_i$, $G\sbe\bigcup_{j=1}^l G_j$ and $\cl_{\mathbb{I}\ap}(f(F_i))\cap \cl_{\mathbb{I}\ap}(f(G_j))=\ems$ for all $i=1,\ldots,k$ and all $j=1,\ldots,l$). This shows that the relations $\s$ and $\ts$ agree on $\BBBB$ (or, equivalently, on $\wt{\BBBB}$).

 Finally, as in every LCA, for every $F,G\in RC(X)$, we have that $F\s G$ iff (there exist $F\ap,G\ap\in\BBBB$ such that $F\ap\sbe F$, $G\ap\sbe G$ and $F\ap\s G\ap$).
 Therefore, the relations $\s$ and $\ts$ agree on $RC(X)$ (or, equivalently, on $\ti{A}$).
 \sqs

 \begin{theorem}\label{whrn}
For every $n\in{\NNNN^+}$, the CLCA $(RC(\mathbb{R}^n),\rho_{\mathbb{R}^n}, CR(\mathbb{R}^n))$ $(=\Psi^t(\mathbb{R}^n))$ is LCA-isomorphic to the $\DHLC$-sum $(\ti{A}_n,\ts_n,\wt{\BBBB}_n)$ of \ $n$ copies of the CLCA $(\ti{A},\ts,\wt{\BBBB})$, constructed in \ref{constr};
thus, the CLCA $(\ti{A}_n,\ts_n,\wt{\BBBB}_n)$ completely determines the Euclidean space $\mathbb{R}^n$ with its natural topology. For every $n\in\NNNN^+$, the Boolean algebras $\ti{A}_n$ and $\ti{A}$ are isomorphic.
\end{theorem}

\doc Since $\mathbb{J}^n$ is homeomorphic to $\mathbb{J}$ and is dense in $\mathbb{R}^n$, we get that $RC(\mathbb{R}^n)$ is isomorphic to $RC(\mathbb{J})$, and thus, to $\ti{A}$ (see \ref{constr} and the proof of Theorem \ref{whr}).  Now all follows from Theorems \ref{whr} and \ref{lcasum}. \sqs

We will now present the description of the CLCA $(RC(\mathbb{R}),\rho_{\mathbb{R}}, CR(\mathbb{R}))$ in two new forms; the notation used in them permits to obtain a more compact form of the definitions of the corresponding relations. As we have already mentioned, $RC(\mathbb{R})$ is isomorphic to $RC(\mathbb{J})$, i.e. to
$RC(\mathbb{Z}_0^{{\NNNN^+}})$  or, equivalently, to $RC(\omega^\omega)$. The last algebra, which is one of the collapsing algebras $RC(k^\omega)$ (where $k$ is an infinite cardinal equipped with the discrete topology), has many abstract descriptions. The one, which is the most appropriate for our purposes, is the following: a complete Boolean algebra $C$ is isomorphic to the Boolean algebra $RC(k^\omega)$ iff it has a dense subset isomorphic to $T^*$, for the normal tree $T=\bigcup\{k^n\st n\in{\NNNN^+}\}$ (here $T^*$ is the tree $T$ with the opposite partial order and $k^n\cap k^m=\ems$ for $n\neq m$) (see, e.g., \cite[14.16(a),(b)]{kop89}).
(Recall that a partially ordered set $(T,\le_T)$ is called a {\em tree}\/ if for every $t\in T$, the set $pred(t)$ is well-ordered by $\le_T$.)
 This shows that $RC(k^\omega)$ is isomorphic to the Boolean algebra $RC(T^*)$, where the ordered set $T^*$ is endowed with the {\em left topology}, i.e. that one generated by the base $\{L_{T^*}(t)\st t\in T \}$ (here $L_{T^*}(t)=\{t\ap\in T\st t\ap\le_{T^*} t\}=\{t\ap\in T\st t\le_{T} t\ap\}$, for every $t\in T$) (see, e.g., \cite[4.11-4.16]{kop89} and \cite[1.7.2]{E}).

Let us add some details and introduce some notation.

\begin{nota}\label{treed}
\rm
For any $n\in{\NNNN^+}$, we set $$\underline{n}=\{1,\ldots,n\}.$$ We set $$T_0=\bigcup\{\mathbb{Z}_0^n\st n\in{\NNNN^+}\},$$
where $\mathbb{Z}_0^n\cap \mathbb{Z}_0^m=\ems$ for $\neq m$.
Any element $t\in\mathbb{Z}_0^n$ is interpreted, as usual, as a function $t:\underline{n} \lra\mathbb{Z}_0$. Further, we let $\bot\sbe t$ and $\bot\neq t$, for any $t\in T_0$; if $n,n\ap\in{\NNNN^+}$, $t\in \mathbb{Z}_0^n$ and
$t\ap\in\mathbb{Z}_0^{n\ap}$, then we set $t\sbe t\ap$ iff $t\ap$ is an extension of $t$, i.e. iff $n\le n\ap$ and $t(i)=t\ap(i)$ for any $i\in\underline{n}$. Then the ordered set $(T_0\cup\{\bot\},\sbe)$ is a normal tree of height $\omega$ with $\mathbb{Z}_0^n$ as its $n$th level (it will be denoted by $L_n$). We also put, for any $t,t\ap\in T_0\cup\{\bot\}$, $$t\le t\ap\Leftrightarrow t\ap\sbe t.$$  We set $$T_0^*=(T_0\cup\{\bot\},\le).$$

 Let $T_0^*$ be endowed with its left topology (i.e. let $(T_0\cup\{\bot\},\sbe)$ be equipped with its right topology (which is defined analogously to the left topology (see \cite[1.7.2]{E}))). Further, for any $t\in T_0\cup\{\bot\}$, put $$c_t=\{t\ap\in T_0\st t\mbox{ and } t\ap\mbox{ are $T_0^*$-compatible}\}.$$ (Recall that two elements $x$ and $y$ of a partially ordered set $(M,\preceq)$ are {\em compatible}\/ if there is some $z\in M$ such that $z\preceq x$ and $z\preceq y$.) Then, as it is well known (see, e.g., \cite[4.13,4.16,the formula for $\cl(u_p)$ in the proof of 4.16]{kop89}), the embedding $e$ of the partially ordered set $T_0^*$ into the Boolean algebra $RC(T_0^*)$ is given by the formula
$$e(t)=c_t,\ \fa t\in T_0\cup\{\bot\}.$$  (Note that the map $e$ is an embedding because $T_0^*$ is a {\em separative partial order} (see, e.g., \cite[4.15,4.16,p.226]{kop89}).) Also, let us recall  that the left topology on  $T_0\cup\{\bot\}$ induced by the ordered set $T_0^*$ is an {\em Alexandroff topology}, i.e. the union of arbitrarily many closed sets is a closed set (see, e.g., \cite[1.7.2]{E}). Thus, the (finite or infinite) joins $\bv\{F_j\st j\in J\}$ in $RC(T_0^*)$ are just the unions $\bigcup\{F_j\st j\in J\}$.

Finally, for every $n\in{\NNNN^+}\stm\{1\}$ and every $t\in L_n$ (i.e. $t:\underline{n}\lra\mathbb{Z}_0$), define
%, for $n> 1$,
\begin{equation}\label{tlam}
t_\l:\underline{n}\lra\mathbb{Z}_0 \mbox{  by the formulas }(t_\l)_{|\ \underline{n-1}}=t_{|\ \underline{n-1}}\mbox{ and } t_\l(n)=(t(n))^-;
\end{equation}
 let, for $t\in L_1$, $t_\l:\underline{1}\lra\mathbb{Z}_0$ be defined by $t_\l(1)=(t(1))^-$.
\end{nota}

\begin{nist}\label{remtreed}
\rm
As we have already mentioned, the Boolean algebra $RC(\mathbb{Z}_0^{{\NNNN^+}})$ is isomorphic to the Boolean algebra $RC(T_0^*)$ (see, e.g., \cite[14.16(a),(b),4.11-4.16]{kop89}). We will recall the proof of this fact since we will use it later. For every $t\in T_0$, set
\begin{equation}\label{at}
a_t=\{x\in \mathbb{Z}_0^{{\NNNN^+}}\st t\sbe x\}.
\end{equation}
Note that if $t:\underline{n} \lra\mathbb{Z}_0$, where $n\in{\NNNN^+}$, then
\begin{equation}\label{att}
a_t=\bigcap_{i=1}^n\pi_i\inv(t(i))
\end{equation}
and thus $a_t$ is a clopen subset of $\mathbb{Z}_0^{{\NNNN^+}}$. Set
\begin{equation}\label{sat}
S=\{a_t\st t\in T_0\}\cup \mathbb{Z}_0^{{\NNNN^+}}.
\end{equation}
 Then $S\sbe CO(\mathbb{Z}_0^{{\NNNN^+}})\sbe RC(\mathbb{Z}_0^{{\NNNN^+}})$.
 Now it is easy to see that the set $S$ is dense in $RC(\mathbb{Z}_0^{{\NNNN^+}})$ and isomorphic to $T_0^*$ (indeed, the map
 \begin{equation}\label{atmap}
 s:T_0^*\lra S, \mbox{ where } s(\bot)=\mathbb{Z}_0^{{\NNNN^+}} \mbox{ and } s(t)=a_t, \fa t\in T_0
 \end{equation}
 is an isomorphism). Therefore, $RC(\mathbb{Z}_0^{{\NNNN^+}})$ is isomorphic to the Boolean algebra $RC(T_0^*)$.
\end{nist}

We will now equip the Boolean algebra $RC(T_0^*)$ defined above with an LCA-structure $( RC(T_0^*),\t,\mathbb{B}_T)$ and will prove that the obtained CLCA is LCA-isomorphic to the CLCA $(RC(\mathbb{R}),\rho_{\mathbb{R}}, CR(\mathbb{R}))$. Recall that two elements $x$ and $y$ of a partially ordered set $(M,\preccurlyeq)$ are {\em comparable}\/ if $x\preccurlyeq y$ or $y\preccurlyeq x$.

\begin{nist}\label{rctost}
\rm
{\bf The construction of the triple} $(RC(T_0^*),\t,\mathbb{B}_T)$.

For every $k,n\in{\NNNN^+}$ and for every $t\in L_k$ (recall that $L_k=\mathbb{Z}_0^k$), set
%$$d_{tn}=\bigcup\{c_{t\ap}\st  (t\ap\in L_{k+1})\mbox{ and }[(t_\l\sbe t\ap, t\ap(k+1)>n)\mbox{ or }(t\sbe t\ap,t\ap(k+1)<-n)]\}.$$
$$d_{tn}=\bigcup\{c_{t\ap}\st  (t\ap\in L_{k+1})\mbox{ \& }[(t_\l\sbe t\ap\ \& \ t\ap(k+1)>n)\mbox{ or }(t\sbe t\ap\ \& \ t\ap(k+1)<-n)]\}.$$
Note that the fact that the left topology on $T_0^*$ is an Alexandroff topology implies that
\begin{equation}\label{dtn}
d_{tn}=
\end{equation}
$$\bv\{c_{t\ap}\st  (t\ap\in L_{k+1})\mbox{ \& }[(t_\l\sbe t\ap\ \mbox{ and } \ t\ap(k+1)>n)\mbox{ or }(t\sbe t\ap\ \mbox{ and } \ t\ap(k+1)<-n)]\}.$$
Let
\begin{equation}\label{C0C1}
C_0=\{c_t\st t\in T_0\}\mbox{ and }C_1=\{d_{tn}\st t\in T_0, \ n\in{\NNNN^+}\}.
\end{equation}

Denote by $\mathbb{B}_{T_0}$ the ideal of $RC(T_0^*)$ generated by $C_0\cup C_1$.

For every $k,k\ap,n,n\ap\in{\NNNN^+}$ and every $t\in L_k$, $t\ap\in L_{k\ap}$, set
%$$
\begin{equation}\label{bbap1}
c_t\t c_{t\ap}\Leftrightarrow
\left\{
\begin{array}{ll}
t=t\ap\mbox{ or }t=t\ap_\l \mbox{ or }t\ap=t_\l, & \mbox{if $k=k\ap$}\\
t\mbox{ and } t\ap\mbox{ are comparable}, & \mbox{if $k\neq k\ap$,}
\end{array}
\right.
\end{equation}
%$$
and
 %$$
 \begin{equation}\label{ububap1}
 d_{tn}\t d_{t\ap n\ap}\Leftrightarrow
 \end{equation}
 $$ \left\{
 \begin{array}{lllll}
  (
 t\ap\sbe t
 \mbox{ and }
%\ \& \
t(k\ap +1)< -n\ap )
 \mbox{ or }(t_\l\ap\sbe t
 \mbox{ and }
%\ \& \
t(k\ap +1) > n\ap ),
 & \mbox{if }k>k\ap +1\\
  (
 t\ap\sbe t
 \mbox{ and }
%\ \& \
t(k)\le -n\ap )
 \mbox{ or }(t_\l\ap\sbe t
 \mbox{ and }
%\ \& \
t(k) > n\ap),
  & \mbox{if }k=k\ap +1\\
 t=t\ap, & \mbox{if $k=k\ap$}\\
  (
 t\sbe t\ap
 \mbox{ and }
 %\ \& \
 t\ap(k\ap)\le -n)
 \mbox{ or }
 (t_\l\sbe t\ap
 \mbox{ and }
 %\ \& \
 t\ap(k\ap) > n),
 & \mbox{if $k=k\ap -1$}\\
 (
 t\sbe t\ap
 \mbox{ and }
 %\ \& \
 t\ap(k+1)< -n)
 \mbox{ or }(t_\l\sbe t\ap
 \mbox{ and }
% \ \& \
t\ap(k+1) > n),
 & \mbox{if $k< k\ap -1$;}
\end{array}\right.
%\end{equation}
$$
and also
%$$
\begin{equation}\label{ubrbap1}
d_{tn}\t c_{t\ap}\Leftrightarrow c_{t\ap}\t d_{tn} \Leftrightarrow
\end{equation}
 $$   \left\{
 \begin{array}{llll}
 t\ap\sbe t, & \mbox{if $k\ap < k$}\\
 t\ap = t\mbox{ or } t\ap=t_\l, & \mbox{if $k\ap=k$}\\
 (t_\l\sbe t\ap\mbox{ and  $t\ap(k\ap)\ge n$) or } (t\sbe t\ap\mbox{ and } t\ap(k\ap)\le -n), & \mbox{if $k\ap=k+1$}\\
 (t_\l\sbe t\ap\mbox{ \& $t\ap(k+1)> n$) or } (t\sbe t\ap \mbox{ \& }t\ap(k+1)< -n), & \mbox{if $k\ap > k+1$.}
\end{array}
\right.
$$
Further, for every two elements $c$ and $d$ of $\mathbb{B}_{T_0}$, set
%$$
\begin{equation}\label{csd1}
c(-\t)d\Leftrightarrow (\ex k,l\in\mathbb{N}^+ \mbox{ and }\ex c_1,\ldots,c_k,d_1,\ldots, d_l\in C_0\cup C_1
\mbox{ such that }
\end{equation}
%$$
$$c\sbe\bigcup_{i=1}^k c_i,\ d\sbe\bigcup_{j=1}^l d_j\mbox{ and } c_i(-\t)d_j,\ \fa i=1,\ldots, k\mbox{ and }\fa j=1,\ldots, l).$$
Finally, for every two elements $a$ and $b$ of $RC(T_0^*)$, set
%$$
\begin{equation}\label{asb1}
a\t b\Leftrightarrow (\ex c,d\in\mathbb{B}_{T_0}\mbox{ such that } c\sbe a,\ d\sbe b\mbox{ and } c\t d).
\end{equation}
%$$
\end{nist}

\begin{theorem}\label{whrct}
The triple $(RC(T_0^*),\t,\mathbb{B}_{T_0})$, constructed in \ref{rctost}, is a CLCA; it is LCA-isomorphic to the complete local contact algebra $(RC(\mathbb{R}),\rho_\mathbb{R},CR(\mathbb{R}))$.
Thus, the triple $(RC(T_0^*),\t,\mathbb{B}_{T_0})$ completely determines the real line $\mathbb{R}$ with its natural topology.
\end{theorem}

\doc
In this proof, we will use the notation introduced in \ref{constr}, \ref{treed}, \ref{remtreed} and \ref{rctost}.
As it follows from  \ref{remtreed} and \cite[the proof of 4.14]{kop89}, there is an isomorphism $h:RC(T_0^*)\lra RC(\mathbb{Z}_0^{{\NNNN^+}})$ defined by the formula
$h(c)=\bigvee_{RC(\mathbb{Z}_0^{{\NNNN^+}})}\{a_{t}\st t\in T_0^*, c_t\sbe c\},$ for every $c\in RC(T_0^*).$
Thus, $h(c_t)=a_t=\bigcap_{i=1}^k\pi_i\inv(t(i))$ and $c_t$ corresponds to $\bigwedge_{i=1}^k\p_i(t(i))$ (see \ref{constr}), where $t\in L_k\sbe T_0^*$ (i.e.,  $t:\unl{k}\lra \mathbb{Z}_0$). This implies that $h(C_0)=B_0\ap=\{a_t\st t\in T_0\}$ and $C_0$ corresponds to $B_0=\{\bigwedge_{i=1}^k\p_i(t(i))\st k\in \NNNN^+, t\in L_k\}$ (see (\ref{C0C1}), (\ref{B0ap}), (\ref{B0})). Note that $t_\l$ corresponds to $b_-$ (see (\ref{tlam}) and (\ref{bu})). Since $h$ is a complete homomorphism, we get that $h(d_{tn})=Q_{a_t n}$ and thus $d_{tn}$ corresponds to $q_{a_t n}$,  for every $k,n\in \NNNN^+$ and every $t\in L_k$ (see (\ref{dtn}), (\ref{QFn}), (\ref{qbn})). Then $h(C_1)=B_1\ap$
and hence $C_1$ corresponds to $B_1$ (see (\ref{C0C1}), (\ref{B1ap}), (\ref{B1})). Hence, $h(\mathbb{B}_{T_0})=\mathbb{B}$ and therefore $\mathbb{B}_{T_0}$ corresponds to $\wt{\BBBB}$ (see the line after (\ref{C0C1}), (\ref{BBBB}) and the paragraph after (\ref{B1ap}), the line after (\ref{B1})). Having  all these facts in mind, we obtain easily that the formula
(\ref{bbap1}) follows from the formula
(\ref{bbap}), (\ref{ububap1}) from (\ref{ububap}), (\ref{ubrbap1}) from (\ref{ubrbap}), (\ref{csd1}) from (\ref{csd}) and (\ref{asb1}) from (\ref{asb}).
This completes the proof of our theorem.
%It follows from the proofs of Theorem \ref{whr} and \cite[14.16,14.15]{kop89}.
\sqs

\begin{theorem}\label{whabstr}
A CLCA  $(M,\mu,\mathbb{M})$
  is LCA-isomorphic to the complete local contact algebra $(RC(\mathbb{R}),\rho_\mathbb{R},CR(\mathbb{R}))$ iff there exists an  embedding (between partially ordered sets)
  $\zeta:T_0^*\lra M$ such that the following two conditions are satisfied:\\
  (a) $\zeta(T_0)$ is dense in $M$, and\\
  (b) let $\zeta(t)=z_t$, for every $t\in T_0$, and let the elements $\widetilde{d_{tn}}$ be defined  by the formula (\ref{dtn}) in which $d_{tn}$ is replaced by $\widetilde{d_{tn}}$, and $c_t$ is replaced by $z_t$;
    then the ideal $\mathbb{M}$ is generated by the set $Z=\zeta(T_0)\cup\{\widetilde{d_{tn}}\st t\in T_0, n\in{\NNNN^+}\}$ and  the formulas (\ref{bbap1}), (\ref{ububap1}), (\ref{ubrbap1}), (\ref{csd}), (\ref{asb}) hold with $\t$ and $\ti{\s}$  replaced by $\mu$, $c_t$ by $z_t$, $d_{tn}$ by $\widetilde{d_{tn}}$, $\widetilde{\mathbb{B}}$ by $\mathbb{M}$,  $B_0\cup B_1$ by $Z$, and $\ti{A}$ by $M$.
\end{theorem}

\doc It follows from  Theorem \ref{whrct} and \cite[4.14,14.16]{kop89}. \sqs

 \section{A Whiteheadian-type description of Tychonoff cubes, spheres and tori}

 \begin{theorem}\label{whsn}
 For every $n\in{\NNNN^+}$, the CNCA $(RC(\mathbb{S}^n),\rho_{\mathbb{S}^n})$ $(=\Psi^t(\mathbb{S}^n))$ is CA-isomorphic to the  CNCA $(\ti{A}_n,C_{\ts_n,\wt{\BBBB}_n})$ (see \ref{whrn} for the LCA $(\ti{A}_n,\ts_n,\wt{\BBBB}_n)$, and \ref{Alexprn} for $C_{\ts_n,\wt{\BBBB}_n}$);
thus, the CNCA $(\ti{A}_n,C_{\ts_n,\wt{\BBBB}_n})$ completely determines the $n$-dimensional sphere   $\mathbb{S}^n$ with its natural topology. Note that $\ti{A}_n$ is isomorphic to $\ti{A}$, for every $n\in\NNNN^+$.
\end{theorem}

\doc  As it follows from the proof of  \cite[Theorem 4.8]{VDDB},  if $X$ is a locally compact Hausdorff space then the complete normal contact algebra $(RC(\a X),\rho_{\a X})$ is CA-isomorphic to the
complete normal contact algebra $(RC(X),C_{\rho_X,CR(X)})$. Now, since $\a\mathbb{R}^n$ is homeomorphic to $\mathbb{S}^n$, our result follows from Theorem \ref{whrn}. \sqs

For every cardinal number $\tau$, denote by $\mathbb{T}^\tau$ the space $(\mathbb{S}^1)^\tau$ (for finite $\tau$, this is just the $\tau$-dimensional torus).

 \begin{theorem}\label{whtt}
For every cardinal number $\tau$, the complete normal contact algebra $(RC(\mathbb{T}^\tau),\rho_{\mathbb{T}^\tau})$ $(=\Psi^t(\mathbb{T}^\tau))$ is CA-isomorphic to the $\DHC$-sum of $\tau$ copies of the CNCA $(\ti{A},C_{\ts,\wt{\BBBB}})$ (see Theorem \ref{whsn} for it);
therefore, this $\DHC$-sum completely determines the space $\mathbb{T}^\tau$.
\end{theorem}

\doc Since the CNCA $(RC(\mathbb{S}^1),\rho_{\mathbb{S}^1})$ is CA-isomorphic to the CNCA $(\ti{A},C_{\ts,\wt{\BBBB}})$ (see Theorem \ref{whsn}), our result follows from Theorem \ref{ncasum}. \sqs

Recall that if $A$ is a Boolean algebra and $a\in A$ then the set
$\downarrow (a)=\{b\in A\st b\le a\}$ endowed with the same meets and joins as in $A$
and with complement $b\ap$ defined by the formula $b\ap=b^*\we
a$, for every $b\le a$, is a Boolean algebra; it is denoted by
$A|a$. If $J=\downarrow (a^*)$ then $A|a$ is isomorphic to the
factor algebra $A/J$; the isomorphism $h:A|a\lra A/J$ is the
following: $h(b)=[b]$, for every $b\le a$ (see, e.g., \cite{kop89}).

In \cite{DII}, we proved the following theorem:

\begin{theorem}\label{conregclo}{\rm \cite[Theorem 6.8]{DII}}
Let $X$ be a locally compact Hausdorff space and $F\in RC(X)$. Set
 $B=RC(X)|F$, $\BBBB\ap=\{G\we F\st G\in CR(X)\}$ and let,
for every $a,b\in B$,  $a\eta b$ iff  $a\rho_X b$ (i.e. $a\cap
b\nes$). Then $(B,\eta,\BBBB\ap)$ is LCA-isomorphic to
$\Psi^t(F)$.
\end{theorem}

Using this assertion, we obtain the following  result:
%theorem together with   Theorem \ref{lccont}, (\ref{deltak}) and Theorem \ref{whr} or Theorem \ref{whabstr} leads to the next two results:

\begin{theorem}\label{un1}
Let $(M,\mu,\mathbb{M})$ be a CLCA which is LCA-isomorphic to the CLCA $(RC(\mathbb{R}),\rho_\mathbb{R},CR(\mathbb{R}))$ and $\zeta:T_0^*\lra M$  be the
embedding described in Theorem \ref{whabstr}. Then, for each $t\in T_0$, the CNCA $(M|\zeta(t),\mu\ap)$, where $\mu\ap$ is the restriction of the relation $\mu$ to
$M|\zeta(t)$, is NCA-isomorphic to the CNCA $(RC(\mathbb{I}),\rho_\mathbb{I})$.
\end{theorem}

\doc By (\ref{deltak}), (\ref{att}) and the beginning of the proof of Theorem \ref{whr}, if $t\in T_0$, i.e. $t: \underline{n} \lra\mathbb{Z}_0$ for some $n\in\NNNN^+$,
then the element $\zeta(t)$ coresponds to the element $\Delta_{t(1)\ldots t(n)}$ of $RC(\mathbb{I}\ap)$ (see also the proofs of theorems \ref{whrct} and  \ref{whabstr}). Since $\Delta_{t(1)\ldots t(n)}$ is homeomorphic to $\mathbb{I}$,
our assertion follows from Theorem \ref{conregclo}.
\sqs

The last theorem shows, in particular, that the following assertion holds:

\begin{theorem}\label{un2}
Let $(\ti{A},\ts,\wt{\BBBB})$ be the CLCA described in  \ref{constr}, $m\in{\NNNN^+}$, $n_1,\ldots,n_m\in\mathbb{Z}_0$, $a_j=\{n_j\}$ for $\jom$, $u=\bigwedge_{j=1}^m\p_j(a_j)$ (see \ref{constr} for $\p_j$) and $B=\ti{A}|u$. Then the CNCA $(B,\ts\ap)$, where $\ts\ap$ is the restriction of the relation $\ts$ to $B$, is NCA-isomorphic to the CNCA $(RC(\mathbb{I}),\rho_\mathbb{I})$. In particular, the CNCA $(RC(\mathbb{I}),\rho_\mathbb{I})$ is NCA-isomorphic to the CNCA $(\ti{A}|\p_1(\{1\}),\ts\ap)$.
\end{theorem}

%Now we will construct a CNCA $(\ti{A},\ts\ap)$ and we will show that it is CA-isomorphic to $\Psi^t(\mathbb{I})$.
A direct description of the CNCA $(RC(\mathbb{I}),\rho_\mathbb{I})$ is given below.

\begin{nist}\label{constri}{\bf The construction of $(\ti{A},\ts\ap)$.}
\rm
We will use the notation from \ref{constr}.

We will define a relation $\ts\ap$ on the Boolean algebra $\ti{A}$  constructed in \ref{constr}.

For every $n\in{\NNNN^+}$, set
$$u_n^\uparrow=\p_1( succ(n))\mbox{ and }u_n^\downarrow=\p_1(pred(-n))$$
and let
$$B_2=\{u_n^\uparrow,u_n^\downarrow\st n\in{\NNNN^+}\}.$$
For every $a,b\in B_0\cup B_1\cup B_2$, set
$$a\ts\ap b\Leftrightarrow a\ts b$$
(see \ref{constr} for the definition of the relation $\ts$). For convenience of the reader, we will write down the corresponding formulae.
For every $n,m\in{\NNNN^+}$,
$$u_n^\uparrow\ts\ap u_m^\uparrow,\ \ u_n^\downarrow\ts\ap u_m^\downarrow\mbox{ \ and \ }u_n^\downarrow(-\ts\ap)u_m^\uparrow.$$
Further, for every $n,r\in{\NNNN^+}$ and every $b=\p_1(a_1)\we \ldots\we \p_k(a_k)\in B_0$, where $a_1=\{m\}$,
\begin{equation}\label{buupd}
b\ts\ap u_n^\uparrow\Leftrightarrow
\left\{
\begin{array}{ll}
m\ge n, & \mbox{if $k=1$}\\
m>n, & \mbox{if $k> 1$}
\end{array}
,\right.
\mbox{ \ \ }
b\ts\ap u_n^\downarrow\Leftrightarrow
\left\{
\begin{array}{ll}
m\le -n, & \mbox{if $k=1$}\\
m<-n, & \mbox{if $k> 1$}
\end{array}
\right.
\end{equation}
and
\begin{equation}\label{quupd}
q_{br}\ts\ap u_n^\uparrow\Leftrightarrow m> n, \ \ q_{br}\ts\ap u_n^\downarrow\Leftrightarrow
\left\{
\begin{array}{ll}
m\le -n, & \mbox{if $k=1$}\\
m<-n, & \mbox{if $k> 1.$}
\end{array}
\right.
\end{equation}
Now, for every $c,d\in\ti{A}$, set
\begin{equation}\label{csdr}
c(-\ts\ap)d\Leftrightarrow (\ex k,l\in\mathbb{N}^+ \mbox{ and }\ex c_1,\ldots,c_k,d_1,\ldots, d_l\in B_0\cup B_1\cup B_2
\mbox{ such  }
\end{equation}
%$$
$$\mbox{that }c\le\bv_{i=1}^k c_i,\ d\le\bv_{j=1}^l d_j\mbox{ and } c_i(-\ts\ap)d_j,\ \fa i=1,\ldots, k\mbox{ and }\fa j=1,\ldots, l).$$
\end{nist}

\begin{theorem}\label{whi}
The pair $(\ti{A},\ts\ap)$, constructed in \ref{constri}, is a complete normal contact algebra; it is CA-isomorphic to the CNCA $(RC(\mathbb{I}),\rho_\mathbb{I})$.
Thus, the pair $(\ti{A},\ts\ap)$ completely determines the closed interval\/ $\mathbb{I}$ with its natural topology.
\end{theorem}

\doc The proof of this assertion is analogous to the proof of Theorem \ref{whr}. We will use in it the notation introduced in \ref{whr}, \ref{constr} and \ref{constri}.

Clearly, $RC(\mathbb{R})$ is isomorphic to $RC(\mathbb{I})$ (by Lemma \ref{isombool}). Thus, $RC(\mathbb{I})$ is isomorphic to $RC(X)$, where $X=\mathbb{Z}_0^{{\NNNN^+}}$ (see the proof of Theorem \ref{whr}).
We will now construct an  NCA $(RC(X),\s\ap)$   CA-isomorphic to $(RC(\mathbb{I}),\rho_{\mathbb{I}})$. Then, identifying $RC(X)$ with $\ti{A}$, we will
  show that $\s\ap=\ts\ap$.

 For every two elements $M$ and $N$ of $RC(\mathbb{J}_2)$, set $M\rho_1 N\Leftrightarrow
\cl_{\mathbb{I}}(M)\cap\cl_{\mathbb{I}}(N)\nes$. Then, using Lemma \ref{isombool}, we get that the pair $(RC(\mathbb{J}_2),\rho_1)$ is CA-isomorphic to the NCA $(RC(\mathbb{I}),\rho_{\mathbb{I}})$.
Now, for every two elements $F,G\in RC(X)$, we set
%$$
\begin{equation}\label{sigrhoi}
F\s\ap G\Leftrightarrow f(F)\rho_1 f(G),
\end{equation}
%$$
where $f:X\lra\mathbb{J}_2$ is the homeomorphism constructed in the proof of Theorem \ref{whr}.
 Obviously, $(RC(X),\s\ap)$  is CA-isomorphic to $(RC(\mathbb{I}),\rho_{\mathbb{I}})$. In the rest of this proof, we will show that the definition of $\s\ap$ given above agrees with the  definition of  $\ti{\s}\ap$ given in \ref{constri}.

 Using the proof of Proposition \ref{rcsemi},
 it is easy to see that the set $$B_2\ap=\{\pi_1\inv( succ(n)),\pi_1\inv(pred(-n))\st n\in{\NNNN^+}\}$$ corresponds to the set $B_2$ introduced in \ref{constri}.
 Now, the formula (\ref{deltak}) implies that, for every $n\in{\NNNN^+}$,
 \begin{equation}\label{ui}
 \ \ \cl_\mathbb{I}(f(\pi_1\inv( succ(n))))=[1-\frac{1}{2^{n+1}},1]\mbox{ and }\cl_\mathbb{I}(f(\pi_1\inv(pred(-n))))=[0,\frac{1}{2^{n+1}}].
 \end{equation}
Thus, for every $m,n\in{\NNNN^+}$, $\cl_\mathbb{I}(f(\pi_1\inv( succ(n))))\cap\cl_\mathbb{I}(f(\pi_1\inv(pred(-m))))=\ems$. Also,
for every $m,n\in{\NNNN^+}$, we have that $f(\pi_1\inv( succ(n)))\cap f(\pi_1\inv( succ(m)))\nes$ and  $f(\pi_1\inv(pred(-n)))\cap f(\pi_1\inv(pred(-m)))\nes.$
Having in mind these formulae and the fact that $\cl_\mathbb{I}(f(F))=\cl_{\mathbb{I}\ap}(f(F))$, for every $F\in B_0\ap\cup B_1\ap$ (see the proof of Theorem \ref{whr} for the notation), we get that $G\s H\Leftrightarrow G\s\ap H$, for every $G,H\in B_0\ap\cup B_1\ap\cup B_2\ap$. This shows that
 $a\ts\ap b\Leftrightarrow a\ts b$, for every $a,b\in B_0\cup B_1\cup B_2$. Hence, the definitions of $\s\ap$ and $\ts\ap$ agree on $B_0\ap\cup B_1\ap\cup B_2\ap$ (or, equivalently, on $B_0\cup B_1\cup B_2$).

Further, using (\ref{ui}), we get that the family $\BB_1=\BB\cup\{\int_\mathbb{I}(\cl_\mathbb{I}(f(F)))\st F\in B_2\ap\}$ (see the proof of Theorem \ref{whr} for the notation and for the fact that $\BB$ is a base of $\mathbb{I}\ap$) is a base of $\mathbb{I}$. Thus, by the regularity of $\mathbb{I}$, every two disjoint elements of $RC(\mathbb{I})$ can be separated by the finite unions of the elements of the family $\{\cl_\mathbb{I}(f(F))\st F\in B_0\ap\cup B_1\ap\cup B_2\ap\}$. This implies that the definitions of $\s\ap$ and $\ts\ap$ agree on $RC(X)$ (or, equivalently, on $\ti{A}$). \sqs

\begin{theorem}\label{whit}
For every cardinal number $\tau$, the complete normal contact algebra $(RC(\mathbb{I}^\tau),\rho_{\mathbb{I}^\tau})$ $(=\Psi^t(\mathbb{I}^\tau))$ is CA-isomorphic to the $\DHC$-sum of $\tau$ copies of the CNCA $(\ti{A},\ts\ap)$ (see Theorem \ref{whi} for it);
therefore, this $\DHC$-sum completely determines the space $\mathbb{I}^\tau$.
\end{theorem}

\doc  It follows from Theorems \ref{whi} and \ref{ncasum}. \sqs

%--------------------------------------------------------------------

\baselineskip = 0.75\normalbaselineskip

\end{document}